\documentclass[12pt,reqno]{amsart}
\usepackage[top=2cm,bottom=2cm,right=2.5cm,left=2.5cm]{geometry}
\usepackage{amssymb}
\usepackage{amsmath, amsthm, amscd, amsfonts, amssymb, graphicx, color}
\usepackage[bookmarksnumbered, colorlinks, plainpages]{hyperref}

\usepackage{hyperref}

\newtheorem{theorem}{Theorem}[section]
\newtheorem{lemma}[theorem]{Lemma}
\newtheorem{remark}[theorem]{Remark}
\newtheorem{proposition}[theorem]{Proposition}
\newtheorem{corollary}[theorem]{Corollary}

\theoremstyle{definition}

\newtheorem{example}[theorem]{Example}

\begin{document}
	
	\setcounter{page}{1}
	
	\title[Admissibility for Multi-window Gabor systems in $\ell^2(\mathbb{S},\mathbb{C}^R)$]{ Admissibility of Multi-window Gabor Systems in Periodically Supported $\ell^2$-spaces with Vector-valued Sequences
	}
	
	\author[N. Khachiaa]{Najib Khachiaa$^*$}
	
	\address{Department of Mathematics Faculty of Sciences, University of Ibn Tofail, B.P. 133, Kenitra, Morocco}
	\email{\textcolor[rgb]{0.00,0.00,0.84}{khachiaa.najib@uit.ac.ma}}
	\date{
		\newline \indent $^{*}$ Corresponding author}

	\begin{abstract}
		In this paper, \( L, M, N, R \) are positive integers, and \( \mathbb{S} \) is an \( N \)-periodic subset of \( \mathbb{Z} \). The space \( \ell^2(\mathbb{S}, \mathbb{C}^R) \) denotes the Hilbert space of vector-valued square-summable sequences over \( \mathbb{S} \), with values in the complex Euclidean space \( \mathbb{C}^R \).
		
		We consider the (multi-window) Gabor system \( \mathcal{G}(g, L, M, N, R) \), generated by applying translations with parameter \( nN \), \( n \in \mathbb{Z} \), and modulations with parameter \( \frac{m}{M} \), \( m \in \mathbb{N}_M \), to a collection of sequences \( g = \{g_l\}_{l \in \mathbb{N}_L} \subset \ell^2(\mathbb{S}, \mathbb{C}^R) \).
		
		Using the vector-valued Zak transform, we characterize the class of sequences \( g \), called windows, that generate a complete Gabor system or a Gabor frame in \( \ell^2(\mathbb{S}, \mathbb{C}^R) \). Furthermore, we provide admissibility conditions under which the periodic set \( \mathbb{S} \) supports a complete Gabor system, a Parseval Gabor frame, or an orthonormal Gabor basis, expressed in terms of the parameters \( L \), \( M \), \( N \), and \( R \).
	\end{abstract}

		\maketitle
	\baselineskip=12.4pt

	\section{Introduction and preliminaries}
	The study of Gabor frames has become a cornerstone of time-frequency analysis, with numerous applications across various scientific and engineering fields. Traditionally, Gabor frames are constructed using a single window, enabling the decomposition of signals into localized time-frequency components. However, the limitations of single-window Gabor frames, particularly in handling signals with diverse characteristics, have led to the development of multi-window Gabor frames. These frames provide greater versatility by employing multiple windows, thereby offering enhanced flexibility in capturing the intricate structures of complex signals. Due to their potential applications in multiplexing techniques, such as Time Division Multiple Access (TDMA) and Frequency Division Multiple Access (FDMA), Vector-valued frames have garnered interest among mathematicians and engineers. In this paper, we focus on multi-window vector-valued Gabor systems.
	
	A sequence $\{f_i\}_{i\in \mathcal{I}}$, where $\mathcal{I}$ is a countable set, in a separable Hilbert space $\mathcal{U}$ is said to be frame if there exist $0< A\leq B<\infty$ (called frame bounds) such that for all $f\in \mathcal{U}$, $$A\|f\|^2\leq \displaystyle{\sum_{i\in \mathcal{I}}\vert \langle f,f_i \rangle\vert^2}\leq B\|f\|^2.$$ If only the upper inequality holds,  $\{f_i\}_{i\in \mathcal{I}}$ is called  a Bessel sequence with Bessel bound $B$. If $A=B$, the sequence is called a $A$-tight frame and if $A=B=1$, it is called a Parseval frame for $\mathcal{U}$. Every frame  $\{f_i\}_{i\in \mathcal{I}}$ for $\mathcal{U}$ defines a positive, self-adjoint, and invertible operator $S$, said frame operator, and defined for all $f\in \mathcal{U}$ by: $S(f):=\displaystyle{\sum_{i\in I}\langle f,f_i\rangle f_i}$. Note that this correspondence is surjective. One of the main strengths of frame theory lies in the redundant and stable decomposition of any signal (i.e., any element of the given Hilbert space) as a countable linear combination of the frame elements. In fact, if $\{f_i\}_{i\in \mathcal{I}}$ is a frame for $\mathcal{U}$ with frame operator $S$, then, for any $f\in \mathcal{U}$, we have the following decomposition: 
	$$f=\sum_{i\in \mathcal{I}}\langle f,S^{-1}f_i\rangle f_i.$$
	For more details on general frame theory, the reader can refer to \cite{11} and  \cite{1}.\\
	Denote by $\mathbb{N}$ the set of positive integers, i.e. $\mathbb{N}:=\{1,2,3,...\}$ and for a given $K\in \mathbb{N}$, write $\mathbb{N}_{K}:=\{0,1,...,K-1\}$. Let $N,M,L\in \mathbb{N}$ and $p,q\in \mathbb{N}$ such that pgcd$(p,q)=1$ and $\displaystyle{\frac{N}{M}=\frac{p}{q}}$. A nonempty subset $\mathbb{S}$ of $\mathbb{Z}$ is said to be $N\mathbb{Z}$-periodic set if for all $j\in \mathbb{S}$ and for all $n\in \mathbb{Z}$, $j+nN \in \mathbb{S}$. For $K\in \mathbb{N}$, write $\mathbb{S}_K:=\mathbb{S}\cap \mathbb{N}_K$. We denote by $\ell^2(\mathbb{S})$ the closed subspace of the Hilbert space $\ell^2(\mathbb{Z})$  of square-summable sequences indexed by the integers, defined by, $$\ell^2(\mathbb{S}):=\{f\in \ell^2(\mathbb{Z}):\, f(j)=0 \text{ if } j\notin  \mathbb{S} \}.$$ 
	Define the modulation operator $E_{\frac{m}{M}}$ with $m\in \mathbb{Z}$ and the translation operator $T_{nN}$ with $n\in \mathbb{Z}$ for $f\in \ell^2(\mathbb{S})$ by: $$E_{\frac{m}{M}}f(.):=e^{2\pi i \frac{m}{M}.} f(.), \;\; T_{nN}f(.):=f(.-nN).$$ 
	The modulation and translation  operators are unitary operators of $\ell^2(\mathbb{S})$.\\ For $g:=\{g_l\}_{l\in \mathbb{N}_L} \subset \ell^2(\mathbb{S})$, the associated Gabor system is given by,
	$$\mathcal{G}(g,L,M,N):=\{E_{\frac{m}{M}}T_{nN}g_l\}_{m\in \mathbb{N}_M,n\in \mathbb{Z}, l\in \mathbb{N}_L}.$$ 
	Let $R\in \mathbb{N}$ and let $\{H_r\}_{r=1,\ldots,R}$ be a sequence of Hilbert spaces on $\mathbb{C}$,  we denote by $\displaystyle{\bigoplus_{r=1}^R H_r}$ their direct sum space equipped with the inner product $\langle f,\tilde{f}\rangle:=\displaystyle{\sum_{r=1}^R\langle f_r,\tilde{f}_r\rangle_{H_r}}$ for $f=(f_1,\ldots,f_R)$ and $\tilde{f}=(\tilde{f}_1,\ldots,\tilde{f}_R)\in \displaystyle{\bigoplus_{r=1}^R H_r}$. In particular, when $H_r=\ell^2(\mathbb{S})$ for all $r=1,\ldots, R$, $\displaystyle{\bigoplus_{r=1}^R H_r}$ is the Hilbert space $\ell^2(\mathbb{S},\mathbb{C}^R)$, and when $H_r=L^2(E)$ for all $r=1,\ldots,R$,  $\displaystyle{\bigoplus_{r=1}^R H_r}$ is the Hilbert space $L^2(E,\mathbb{C}^R)$. For  details  frame theory on the direct sum of Hilbert spaces, the reader can refer to \cite{2}. In this paper, we are interested in Gabor frames for  $\ell^2(\mathbb{S},\mathbb{C}^R)$.\\
	
	Define the modulation operator $E_{\frac{m}{M}}$ with $m\in \mathbb{Z}$ and the translation operator $T_{nN}$ with $n\in \mathbb{Z}$ for $f\in \ell^2(\mathbb{S}, \mathbb{C}^R)$ by: $$E_{\frac{m}{M}}f(.):=(E_{\frac{m}{M}}f_1(.),E_{\frac{m}{M}}f_2(.),\ldots,E_{\frac{m}{M}}f_R(.)\,)$$
	and $$ T_{nN}f(.):=(T_{nN}f_1(.),T_{nN}f_2(.)\ldots,T_{nN}f_R(.)\,).$$ 
	The modulation and translation  operators are unitary operators of $\ell^2(\mathbb{S},\mathbb{C}^R)$.\\ For $g:=\{g_l\}_{l\in \mathbb{N}_L} \subset \ell^2(\mathbb{S},\mathbb{C}^R)$, the associated  Gabor system is given by,
	$$\mathcal{G}(g,L,M,N,R):=\{E_{\frac{m}{M}}T_{nN}g_l\}_{m\in \mathbb{N}_M,n\in \mathbb{Z}, l\in \mathbb{N}_L}.$$ 
	If $\mathcal{G}(g,L,M,N,R)$ is a  Gabor frame and $S\in B(\ell^2(\mathbb{S},\mathbb{C}^R)\,)$ is its frame  operator, then it is easy to verify that $S^{-1}(\mathcal{G}(g,L,M,N,R))=\mathcal{G}(S^{-1}g,L,M,N,R))$ and $S^{\frac{-1}{2}}(\mathcal{G}(g,L,M,N,R))=\mathcal{G}(S^{\frac{-1}{2}}g,L,M,N,R))$ which is a Parseval  Gabor frame.
	Let $K\in \mathbb{N}$, the discrete Zak tansform $z_K$ of $h\in \ell^2(\mathbb{Z})$ for $j\in \mathbb{Z}$ and a.e $\theta\in \mathbb{R}$ is defined by, 	
	$$z_Kh(j,\theta):=\displaystyle{\sum_{k\in \mathbb{Z}}h(j+kK)e^{2\pi i k \theta}}.$$
	For $f\in \ell^2(\mathbb{Z},\mathbb{C}^R)$, the vector-valued discrete Zak transform $z_Kf$ for $f$ is defined by $$
	z_Kf(j,\theta):=(z_kf_1(j,\theta),z_Kf_2(j,\theta),\ldots,z_Kf_R(j,\theta)\,).$$
	We can verify easily that the vector-valued Zak transform $z_Kf$ is  quasi-periodic, i.e. for all $j,k,l\in \mathbb{Z}$ and $\theta\in \mathbb{R}$, we have: $$z_Kf(j+kK,\theta+l)=e^{-2\pi ik\theta}z_Kf(j,\theta).$$
	Then, the vector-valued discrete Zak transform is  completely determined by its values for $j\in \mathbb{N}_K$ and $\theta \in [0,1[$.\\
	
	This paper is organized as follows. In section 2,  we will present some auxiliary lemmas to be used in the following sections. In section 3, we characterize which $g\in \ell^2(\mathbb{S},\mathbb{C}^R)$ generates a complete  Gabor system and a  Gabor frame $\mathcal{G}(g,L,M,N,R)$ for $\ell^2(\mathbb{S},\mathbb{C}^R)$ using the vector-valued Zak transform. In section 4, we provide an admissibility characterization for complete Gabor systems,  Gabor( Parseval) frames and  Gabor (orthonormal) bases $\mathcal{G}(g,L,M,N,R)$ on a discrete periodic subset $\mathbb{S}$ of $\mathbb{Z}$, and we conclude with an illustrative example.
	\section{Auxiliary lemmas}	
	In this section, we introduce a series of lemmas and establish the notations that will be used in the following sections. Beyond the notations mentioned in the introduction, let $\mathcal{M}_{s,t}$ denote the collection of all $s\times t$ matrices with elements in $\mathbb{C}$. For $s\in \mathbb{N}$, $I_s$ denotes the identity matrix in $\mathcal{M}_{s,s}$. The notation $p\wedge q=1$ is used to indicate that $p$ and $q$ are coprime. For any matrix $A$, $A^*$ stands for its conjugate transpose, $N(A)$ represents its null space, and $A_{s,t}$ indicates its $(s,t)$ entry. When $A$ is a column vector, we denote its $r$-th entry by $A_r$. Let $L,M,N,R\in \mathbb{N}$ and $p,q\in \mathbb{N}$ such that $p\wedge q=1$ and  that $\displaystyle{\frac{N}{M}=\frac{p}{q}}$. For $j\in \mathbb{Z}$, we denote $\mathcal{K}_j:=\{k\in \mathbb{N}_p:\, j+kM\in \mathbb{S}\}$,  $\mathcal{K}(j):=diag(\chi_{\mathcal{K}_j}(0),\chi_{\mathcal{K}_j}(1),...,\chi_{\mathcal{K}_j}(p-1))$ and $\Lambda_j:=\{k+rp:\; k\in \mathcal{K}_j,\; r\in \mathbb{N}_R\}$. 
	Following this, we provide several definitions and results that will be pertinent throughout the rest of the paper.\\
	
	The following lemma is a vector version of Theorem $2.1$ in \cite{4}.
	\begin{lemma}\label{lem2.1}
		Let $K,R\in \mathbb{N}$, and let $\mathbb{S}$ be a $K\mathbb{Z}$-periodic set in $\mathbb{Z}$. Then the vector-valued Zak transform $z_K$ is a unitary operator from $\ell^2(\mathbb{S},\mathbb{C}^R)$ onto $L^2(\mathbb{S}_K\times [0,1[, \mathbb{C}^R)$, where $$L^2(\mathbb{S}_K\times [0,1[):=\left\{\psi:\mathbb{S}_K\times [0,1[\rightarrow \mathbb{C}:\, \displaystyle{\sum_{j\in \mathbb{S}_K}\int_{0}^1 \vert \psi(j,\theta)\vert^2 \, d\theta}< \infty\right\}.$$
	\end{lemma}
	
	Let $A,B\subset \mathbb{Z}$ and $K\in \mathbb{N}$.  We say that $A$ is $K\mathbb{Z}$-congruent to $B$ if there exists a partition $\{A_k\}_{k\in \mathbb{Z} }$ of $A$ such that $\{A_k+kK\}_{k\in \mathbb{Z}}$ is a partition of $B$.
	\begin{lemma}\label{lem2.2}\cite{4}
		Let $N,M\in \mathbb{N}$ and $p,q\in \mathbb{N}$ such that $p\wedge q=1$ and $\displaystyle{\frac{N}{M}=\frac{p}{q}}.$ Then 
		the set \\$\Delta:=\{j+kM-rN:\; j\in \mathbb{N}_{\frac{M}{q}},\, k\in \mathbb{N}_p,\, r\in \mathbb{N}_q\}$ is $pM$-congruent to $\mathbb{N}_{pM}$.
	\end{lemma}	
	
	Let $R,M,N\in \mathbb{N}$. Write $\displaystyle{\frac{N}{M}}=\displaystyle{\frac{p}{q}}$ such that $p,q\in \mathbb{N}$ and  $p\wedge q=1$. For each $h\in \ell^2(\mathbb{Z})$, we associate a matrix-valued function $Z_h:\mathbb{Z}\times \mathbb{R}\rightarrow \mathcal{M}_{q,p}$ whose entry at the r-th row and the k-th column is defined for $(j,\theta)\in \mathbb{Z}\times \mathbb{R}$ by:
	$$Z_h(j,\theta)_{r,k}=z_{pM}h(j+kM-rN,\theta).$$
	For $f=(f_1,f_2,\ldots,f_R)\in \ell^2(\mathbb{Z}, \mathbb{C}^R)$, we associate a matrix-valued function $Z_f:\mathbb{Z}\times \mathbb{R}\rightarrow\mathcal{M}_{q,pR}$ defined for $(j,\theta)\in \mathbb{Z}\times \mathbb{R}$ by: 
	$$Z_f(j,\theta):=(Z_{f_1}(j,\theta),Z_{f_2}(j,\theta),\ldots,Z_{f_R}(j,\theta)\,).$$
	Combining Lemma $\ref{lem2.2}$ with the quasi-periodicity of the discrete Zak transform, we have the following lemma.
	\begin{lemma}
		Let $f=(f_1,f_2,\ldots,f_R)\in \ell^2(\mathbb{S},\mathbb{C}^R)$ and  let $N,M\in \mathbb{N}$ and $p,q\in \mathbb{N}$ such that $\displaystyle{\frac{N}{M}=\frac{p}{q}}.$ Then 
		$z_{pM}f$ is completeley determined by the matrices $Z_f(j,\theta)$ for $j\in \mathbb{N}_{\frac{M}{q}}$ and $\theta \in [0,1[$. Conversely, a matrix-valued function $Z:\mathbb{N}_{\frac{M}{q}}\times [0,1[\rightarrow \mathcal{M}_{q,pR}$ such that for all $j\in \mathbb{N}_{\frac{M}{q}}$, $Z(j,.)_{r,k}\in L^2([0,1[)$ also determines a unique $f\in \ell^2(\mathbb{Z},\mathbb{C}^R)$ such that for all $j\in \mathbb{N}_{\frac{M}{q}}$, $\theta \in [0,1[$, $Z_f(j,\theta)=Z(j,\theta).$\\
	\end{lemma}
	
	For $g:=\{g_l\}_{l\in \mathbb{N}_L}\subset \ell^2(\mathbb{Z},\mathbb{C}^R)$,  we associate a matrix-valued function $Z_g:\mathbb{Z}\times \mathbb{R}\rightarrow\mathcal{M}_{qL,pR}$ defined for $(j,\theta)\in \mathbb{Z}\times \mathbb{R}$ by:
	$$Z_g(j,\theta)=\begin{pmatrix}
		Z_{g_0}(j,\theta)\\
		Z_{g_1}(j,\theta)\\
		\vdots\\
		Z_{g_{L-1}}(j,\theta)
	\end{pmatrix}=\left(Z_{g_{l,r}}(j,\theta)\right)_{l=0,1,\ldots,L-1,\;r=1,2,\ldots,R}
	$$
	where for all $l\in \mathbb{N}_L$, we write $g_l=(g_{l,1},g_{l,2},\ldots,g_{l,R})$.\\
	
	For $f=(f_1,f_2,\ldots,f_R)\in \ell^2(\mathbb{Z},\mathbb{C}^R)$, we, also, associate a vector-valued function $F:\mathbb{Z}\times \mathbb{R}\rightarrow \mathbb{C}^{pR}$ defined for $(j,\theta)\in \mathbb{Z}\times \mathbb{R}$ by:
	$$F(j,\theta)=\begin{pmatrix}
		F_1(j,\theta)\\
		F_2(j,\theta)\\
		\vdots\\
		F_R(j,\theta)
	\end{pmatrix},
	$$
	where $F_r: \mathbb{Z}\times \mathbb{R}\rightarrow \mathbb{C}^p$ is a vector-valued function defined for $(j,\theta)\in \mathbb{Z}\times \mathbb{R}$ by the column vector: $$(F_r(j,\theta)):=\left(z_{pM}f_r(j+kM,\theta)\right)_{k\in \mathbb{N}_p}.$$
	By Lemma~\ref{lem2.1} and the quasi-periodicity of the discrete Zak transform, the sequence \( f \) and its Zak transform \( F(j, \theta) \), with \( j \in \mathbb{N}_M \) and \( \theta \in [0,1[ \), uniquely determine each other.

	By Lemma 3 in \cite{8}, we obtain the following lemma.
	\begin{lemma}\label{lem2.4}
		Let $M,N\in \mathbb{N}$ and $p,q\in \mathbb{N}$ such that $p\wedge q=1$ and  that $\displaystyle{\frac{N}{M}=\frac{p}{q}}$. Let $h\in \ell^2(\mathbb{Z},\mathbb{C}^R)$, then for all $f\in \ell^2(\mathbb{Z},\mathbb{C}^R)$ and $(m,n,r)\in \mathbb{N}_M\times \mathbb{Z}\times \mathbb{N}_q$, we have:
		$$\langle f, E_{\frac{m}{M}}T_{(nq+r)N}h\rangle =\displaystyle{\sum_{j\in \mathbb{N}_M}\left( \int_0^1 (\overline{Z_h(j,\theta)}F(j,\theta))_re^{-2\pi i n\theta}d\theta\right) e^{-2\pi i \frac{m}{M}j}}.$$
	\end{lemma}
	Let $s,t\in \mathbb{N}$ and  $A\in \mathcal{M}_{s,t}$. For $K\in \mathbb{N}$, we denote by $I_K\otimes A$ the block matrix in $\mathcal{M}_{sK,tK}$ defined by:
	$$\begin{pmatrix}
		A & 0 & \cdots & 0 \\
		0 & A & \cdots & 0 \\
		\vdots & \vdots & \ddots & \vdots \\
		0 & 0 & \cdots & A \\
	\end{pmatrix}$$ 
	
	By  simple computations, we obtain the following lemmas \ref{lem2.5}, \ref{lem2.6}, \ref{lem2.7},  \ref{lem2.9}, \ref{lem2.10}.
	\begin{lemma}\label{lem2.5}
		Let $A\in \mathcal{M}_{s,s}$, $A'\in \mathcal{M}_{t,t}$ and $\{B_{k,k'}\}_{k\in \mathbb{N}_K, k'\in \mathbb{N}_{K'}}\subset \mathcal{M}_{s,t}$. Then:
		$$
		\begin{array}{rcl}
			&&\begin{pmatrix}
				AB_{0,0}A' & AB_{0,1}A' & \cdots & AB_{0,K'-1}A' \\
				AB_{1,0}A'& AB_{1,1}A' & \cdots & AB_{1,K'-1}A' \\
				\vdots & \vdots & \ddots & \vdots \\
				AB_{K-1,0}A' & AB_{K-1,1}A' & \cdots & AB_{K-1,K'-1}A' \\
			\end{pmatrix}\\
			&=&(I_K\otimes A)\begin{pmatrix}
				B_{0,0} & B_{0,1} & \cdots & B_{0,K'-1} \\
				B_{1,0} & B_{1,1} & \cdots & B_{1,K'-1} \\
				\vdots & \vdots & \ddots & \vdots \\
				B_{K-1,0} & B_{K-1,1} & \cdots & B_{K-1,K'-1} \\
			\end{pmatrix}(I_{K'}\otimes A').
		\end{array}$$
	\end{lemma}
	\begin{lemma}\label{lem2.6}
		Let $K\in \mathbb{N}$, $\lambda \in \mathbb{C}$, and $A, B$ be two matrices. Then:
		$$I_K\otimes (\lambda A+B)=\lambda (I_K\otimes A)+(I_K\otimes B).$$
	\end{lemma}
	\begin{lemma}\label{lem2.7}
		Let $K\in \mathbb{N}$ and $A, B$ be two matrices such that $AB$ is well defined. Then:
		$$I_K\otimes (AB)=(I_K\otimes A)(I_K\otimes B).$$
	\end{lemma}
	\begin{remark}\label{lem2.8}
		Let $K,s\in \mathbb{N}$. Then  for all $s\in \mathbb{N}$, $I_K\otimes I_s=I_{sK}$.
	\end{remark}
	\begin{lemma}\label{lem2.9}
		Let $K\in \mathbb{N}$ and  $A$ be a matrix. Then: 
		$$(I_K\otimes A)^*=I_K\otimes A^*.$$
	\end{lemma}
	\begin{lemma}\label{lem2.10}
		Let $K\in \mathbb{N}$ and $A$ be a square matrix. Then: 
		\begin{enumerate}
			\item $A$ is invertible $\Longleftrightarrow$ $I_K\otimes A$ is invertible. Moreover, in this case, we have: 
			$$(I_K\otimes A)^{-1}=I_K\otimes A^{-1}.$$
			\item $A$ is a unitary matrix $\Longleftrightarrow$ $I_K\otimes A$ is a unitary matrix.
		\end{enumerate}
	\end{lemma}
	
	Let $p,q\in \mathbb{N}$ such that $p\wedge q=1$. By Lemma $4$ in \cite{4}, for each $\ell\in \mathbb{Z}$, there exists a unique $(k_{\ell},r_{\ell},m_{\ell})\in \mathbb{N}_p\times \mathbb{N}_q\times \mathbb{Z}$ such that $\ell=k_{\ell}q+(m_{\ell}q-r_{\ell})p$. Define the two following matrix-valued functions as follows ( for a.e $\theta\in [0,1[$):
	$$
	A_{\ell}(\theta):=\left\lbrace
	\begin{array}{rcl}
		&\begin{pmatrix}
			0 & e^{-2\pi i\theta}I_{k_{\ell}}\\
			I_{p-k_{\ell}} & 0
		\end{pmatrix}& \; \text{ if } k_{\ell}\neq 0,\\
		&I_p&\; \text{ if } k_{\ell}=0.
	\end{array}
	\right.$$
	And: 
	$$
	C_{\ell}(\theta):=\left\lbrace
	\begin{array}{rcl}
		&\begin{pmatrix}
			0 &   I_{q-r_{\ell}}\\
			e^{2\pi i\theta}I_{r_{\ell}}  & 0
		\end{pmatrix}& \; \text{ if } r_{\ell}\neq 0,\\
		&I_q&\; \text{ if } k_{\ell}=0.
	\end{array}
	\right.$$
	\begin{lemma}\label{lem2.11}
		The matrices $A_{\ell}(\theta)$ and $C_{\ell}(\theta)$ are unitaries.\\
	\end{lemma}
	\begin{proof}
		By a simple computation, we obtain: 
		$$\begin{array}{rcl}
			A_{\ell}(\theta)A_{\ell}(\theta)^*&=&\begin{pmatrix}
				I_{k_{\ell}} & 0\\
				0 & I_{p-k_{\ell}}
			\end{pmatrix}\\
			&=& I_p.
		\end{array}$$
		And $$\begin{array}{rcl}
			C_{\ell}(\theta)C_{\ell}(\theta)^*&=&\begin{pmatrix}
				I_{q-r_{\ell}} & 0\\
				0 &  I_{r_{\ell}}
			\end{pmatrix}\\
			&=& I_q.
		\end{array}$$
	\end{proof}
	\begin{lemma}\label{lem2.12}
		Let $M,N\in \mathbb{N}$ and $p,q\in \mathbb{N}$ such that $p\wedge q=1$ and  that $\displaystyle{\frac{N}{M}=\frac{p}{q}}$. Let $g:=\{g_l\}_{l\in \mathbb{N}_L}\subset \ell^2(\mathbb{S},\mathbb{C}^R)$. Then, for all $j,\ell\in \mathbb{Z}$ and a.e $\theta \in \mathbb{R}$, we have:
		\begin{equation}
			Z_g(j+\displaystyle{\frac{M}{q}}\ell,\theta)=e^{-2\pi i m_{\ell}\theta}(I_L\otimes C_{\ell}(\theta))Z_g(j,\theta)(I_R\otimes A_{\ell}(\theta)).
		\end{equation}
		Moreover, we have: 
		\begin{equation}
			Z_g(j+\displaystyle{\frac{M}{q}}\ell,\theta)^*Z_g(j+\displaystyle{\frac{M}{q}}\ell,\theta)=(I_R\otimes A_{\ell}(\theta))^*Z_g(j,\theta)^*Z_g(j,\theta)(I_R\otimes A_{\ell}(\theta)).
		\end{equation}
	\end{lemma}
	\begin{proof}
		By Lemma $5$ in \cite{8}, we have for all $l\in \mathbb{N}_L$, $r\in \{1,2,\ldots,R\}$, $j,\ell\in \mathbb{Z}$ and a.e $\theta\in \mathbb{R}$: 
		$$Z_{g_{l,r}}(j+\displaystyle{\frac{M}{q}}\ell,\theta)=e^{-2\pi i m_{\ell}\theta}C_{\ell}(\theta)Z_{g_{l,r}}(j,\theta)A_{\ell}(\theta).$$
		Then, Lemma \ref{lem2.5} implies $(1)$ and the unitarity of $C_{\ell}(\theta)$ together with Lemma \ref{lem2.10} implies $(2)$.
	\end{proof}
	
	\begin{lemma}
		Let $M,N\in \mathbb{N}$ and $p,q\in \mathbb{N}$ such that $p\wedge q=1$ and  that $\displaystyle{\frac{N}{M}=\frac{p}{q}}$. Let $g:=\{g_l\}_{l\in \mathbb{N}_L}\subset \ell^2(\mathbb{Z},\mathbb{C}^R)$. The following statements are equivalent.
		\begin{enumerate}
			\item $\mathcal{G}(g,L,M,N,R)$ is a Bessel sequence in $\ell^2(\mathbb{Z},\mathbb{C}^R)$.
			\item All entries of $Z_g(j,.)$ are in $L^{\infty}([0,1[)$ for each $j\in \mathbb{N}_{\frac{M}{q}}$.
		\end{enumerate}
	\end{lemma}
	
	\begin{proof}
		It is well known that $\mathcal{G}(g,L,M,N,R)$ is a Bessel sequence in $\ell^2(\mathbb{Z},\mathbb{C}^R)$ if and only if for each $1\leq r\leq R$, $\mathcal{G}( \{g_{l,r}\}_{l\in \mathbb{N}_L}\,  ,L,M,N)$ is a Bessel sequence in $\ell^2(\mathbb{Z})$. Since for $1\leq r\leq R$, $\mathcal{G}( \{g_{l,r}\}_{l\in \mathbb{N}_L}\,  ,L,M,N)$ is a Bessel sequence in $\ell^2(\mathbb{Z})$  if and only if $\mathcal{G}(g_{l,r}, M,N)$ is a Bessel sequence in $\ell^2(\mathbb{Z})$ for all $l\in \mathbb{N}_L$, then by Proposition $1$ in \cite{8}, for all $l\in \mathbb{N}_L$ and $1\leq r\leq R$, 
		all entries of $Z_{g_{l,r}}(j,.)$ are in $L^{\infty}([0,1[)$ (for all $j\in \mathbb{N}_{\frac{M}{q}}$). Hence $\mathcal{G}(g,L,M,N,R)$ is a Bessel sequence in $\ell^2(\mathbb{Z},\mathbb{C}^R)$ if and only if all entries of $Z_g(j,.)$ are in $L^{\infty}([0,1[)$ for all $j\in \mathbb{N}_{\frac{M}{q}}$.
	\end{proof}

	\begin{lemma}\label{lem13}
		Let \( g := \{g_l\}_{l \in \mathbb{N}_L} \subset \ell^2(\mathbb{Z}, \mathbb{C}^R) \) and let \( f \in \ell^2(\mathbb{Z}, \mathbb{C}^R) \). Then the following statements are equivalent:
		\begin{enumerate}
			\item \( f \) is orthogonal to the Gabor system \( \mathcal{G}(g, L, M, N, R) \).
			\item \( Z_g(j, \theta)\, \overline{F(j, \theta)} = 0 \) for all \( j \in \mathbb{N}_M \) and almost every  \( \theta \in [0,1[ \).
			\item \( Z_{g_l}(j, \theta)\, \overline{F(j, \theta)} = 0 \) for all \( l \in \mathbb{N}_L \), \( j \in \mathbb{N}_M \), and almost every \( \theta \in [0,1[ \).
		\end{enumerate}
	\end{lemma}

	\begin{proof}
		For $f\in \ell^2(\mathbb{Z},\mathbb{C}^R)$, $f$ is orthogonal to $\mathcal{G}(g,L,M,N,R)$ if and only if $f$ is orthogonal to $\mathcal{G}(g_l,M,N,R)$ for all $l\in \mathbb{N}_L$. Then, by Lemma \ref{lem2.4}, $f\in \ell^2(\mathbb{Z},\mathbb{C}^R)$, $f$ is orthogonal to $\mathcal{G}(g,L,M,N,R)$ if and only if $Z_{g_l}(j,\theta)\overline{F(j,\theta)}=0$ for all $l\in \mathbb{N}_L$, $j\in \mathbb{N}_{M}$ and a.e $\theta \in  [0,1[$, hence $(1)\Longleftrightarrow (3)$. Since $$Z_g(j,\theta)\overline{F(j,\theta)}=\begin{pmatrix}
			Z_{g_0}(j,\theta)\overline{F(j,\theta)}\\
			Z_{g_1}(j,\theta)\overline{F(j,\theta)}\\
			\vdots\\
			Z_{g_{L-1}}(j,\theta)\overline{F(j,\theta)}\\
		\end{pmatrix},$$
		then $(2)\Longleftrightarrow (3).$
	\end{proof}
	\begin{lemma}\label{lem2.14}
		For all $j\in \mathbb{Z}$, $\mathcal{K}(j)$ is an orthogonal projection on $\mathbb{C}^p$. i.e.
		
		$$\mathcal{K}(j)^2=\mathcal{K}(j), \quad 
		\mathcal{K}(j)^*=\mathcal{K}(j).$$
		
	\end{lemma}
	\begin{proof}
		Let $j\in\mathbb{Z}$. We have:  $$\begin{array}{rcl}
			\mathcal{K}(j)^2&=&\text{diag}(\chi_{\mathcal{K}_j}(0)^2,\chi_{\mathcal{K}_j}(1)^2,\ldots,\chi_{\mathcal{K}_j}(p-1)^2)\\
			&=&\text{diag}(\chi_{\mathcal{K}_j}(0),\chi_{\mathcal{K}_j}(1),\ldots, \chi_{\mathcal{K}_j}(p-1))=\mathcal{K}(j).
		\end{array}$$
		And $$\begin{array}{rcl}
			\mathcal{K}(j)^*&=&\text{diag}(\overline{\chi_{\mathcal{K}_j}(0)},\overline{\chi_{\mathcal{K}_j}(1)},\ldots,\overline{\chi_{\mathcal{K}_j}(p-1)})\\
			&=&\text{diag}(\chi_{\mathcal{K}_j}(0),\chi_{\mathcal{K}_j}(1),\ldots, \chi_{\mathcal{K}_j}(p-1))=\mathcal{K}(j).
		\end{array}$$
	\end{proof}
	\begin{lemma}\label{lem2.15}
		Let $j\in \mathbb{Z}$. Then: 
		$$diag\left(\chi_{\Lambda_j}(0),\chi_{\Lambda_j}(1),\ldots,\chi_{\Lambda_j}(pR-1)\,\right)=(I_R\otimes \mathcal{K}(j)\,).$$
	\end{lemma}
	\begin{proof}
		Let $k\in \mathbb{N}_{pR}$, then there exist a unique $(k',r)\in \mathbb{N}_p\times \mathbb{N}_R$ such that $k=k'+pr$. It is clear that $\chi_{\Lambda_j}(k)=\chi_{\mathcal{K}_j}(k')$. Then, Lemma \ref{lem2.5} completes the proof.
	\end{proof}
	\begin{lemma}\label{lem2.16}\cite{6}
		For all $j,\ell\in \mathbb{Z}$ and a.e $\theta\in [0,1[$, we have:
		$$I_R\otimes \mathcal{K}(j+\displaystyle{\frac{M}{q}}\ell)=(I_R\otimes A_{\ell}(\theta)\,)^*(I_R\otimes \mathcal{K}(j)\,)(I_R\otimes A_{\ell}(\theta)\,).$$
	\end{lemma}
	\begin{lemma}\cite{1}\label{lem2.17}
		Let $\{f_i\}_{i\in \mathcal{I}}$, where $\mathcal{I}$ is a countable sequence,  be a Parseval frame for a separable Hilbert space $\mathcal{U}$. Then the following statements are equivalent:
		\begin{enumerate}
			\item $\{f_i\}_{i\in \mathcal{I}}$ is a Riesz basis.
			\item $\{f_i\}_{i\in \mathcal{I}}$ is an orthonormal basis.
			\item For all $i\in \mathcal{I}$, $\|f_i\|=1.$
		\end{enumerate}
	\end{lemma}
	\begin{lemma}\label{lem2.18}
		Let $g:=\{g_l\}_{l\in \mathbb{N}_L}\subset \ell^2(\mathbb{S})$ such that $\mathcal{G}(g,L,M,N)$ is a Parseval frame for $\ell^2(\mathbb{S})$. Then:
		$$\displaystyle{\sum_{l\in \mathbb{N}_L} \|g_l\|^2=\frac{card(\mathbb{S}_N)}{M}}.$$
	\end{lemma}
	\begin{proof}
		By Theorem $3.2$ in \cite{7}, we have for all $j\in \mathbb{S}_N$, $\displaystyle{\sum_{l\in \mathbb{N}_L}\sum_{n\in \mathbb{Z}}\vert g_l(j-nN)\vert^2=\frac{1}{M}}$. Then $\displaystyle{\sum_{j\in \mathbb{S}_N}\sum_{l\in \mathbb{N}_L}\sum_{n\in \mathbb{Z}}\vert g_l(j-nN)\vert^2=\frac{card(\mathbb{S}_N)}{M}}$, thus $\displaystyle{\sum_{l\in \mathbb{N}_L}\|g_l\|^2=\frac{card(\mathbb{S}_N)}{M}}$.
	\end{proof}
	\begin{lemma}\cite{4}\label{lemma6}
		Let $M \in \mathbb{N}$ and $E \subset \mathbb{Z}$. Then the following conditions are equivalent:
		\begin{enumerate}
			\item $\left\{ e^{2\pi im/M} \cdot \chi_E(\cdot) : m \in \mathbb{N}_M \right\}$ is a tight frame for $\ell^2(E)$ with frame bound $M$.
			\item $\left\{ e^{2\pi im/M} \cdot \chi_E(\cdot) : m \in \mathbb{N}_M \right\}$ is complete in $\ell^2(E)$.
			\item $E$ is $M\mathbb{Z}$-congruent to a subset of $\mathbb{N}_M$.
			\item $\displaystyle{\sum_{k \in \mathbb{Z}} \chi_E(\cdot + kM) \leq 1}$ on $\mathbb{Z}$.
		\end{enumerate}
	\end{lemma}

	\section{Characterizations of complete  Gabor systems and Gabor frames for $\ell^2(\mathbb{S},\mathbb{C}^R)$}
	In this section we use all the notations already introduced without introducing them again. Note also that, since no confusion is possible, all the norms that will be used will be noted by the same notation $\|.\|$. Let $L,M,N,R\in \mathbb{N}$ and $p,q\in \mathbb{N}$ such that $p\wedge q=1$ and  $\displaystyle{\frac{N}{M}=\frac{p}{q}}$ and denote for $K\in \mathbb{N}_K$, $\mathbb{S}_K=\mathbb{S}\cap \mathbb{N}_K$. In what follows, we characterize the class of sequences \( g := \{g_l\}_{l \in \mathbb{N}_L} \subset \ell^2(\mathbb{S}, \mathbb{C}^R) \) that generate a complete Gabor system, a Gabor frame, or a Gabor (orthonormal) basis \( \mathcal{G}(g, L, M, N, R) \).

	We, first, present the following proposition.
	\begin{proposition}\label{prop1}
		Let $g:=(g_1,g_2,\ldots,g_R),\; h:=(h_1,h_2,\ldots,h_R)\in \ell^2(\mathbb{Z},\mathbb{C}^R)$ such that $\mathcal{G}(g,L,M,N,R)$ and $\mathcal{G}(h,L,M,N,R)$ are both Bessel sequences in $\ell^2(\mathbb{Z},\mathbb{C}^R
		)$. Then:
		$$\begin{array}{rcl}
			&&\left((z_{pM}S_{h,g}f)(j+kM,\theta)\,\right)_{k\in \mathbb{N}_p}\\
			&=&\displaystyle{M\sum_{l\in \mathbb{N}_L} \left(Z_{g_{l,1}}^t(j,\theta),Z_{g_{l,2}}^t(j,\theta),\ldots,Z_{g_{l,R}}^t(j,\theta)\right)\left(I_R \otimes (\,\overline{Z_{h_l}(j,\theta)}F(j,\theta)\,)\,\right)},
		\end{array}$$
		for $f=(f_1,f_2,\ldots,f_R)\in \ell^2(\mathbb{Z},\mathbb{C}^R)$, $j\in \mathbb{N}_M$ and a.e $\theta \in [0,1[$.\\
	\end{proposition}
	\begin{proof}
		We have: 
		$$\begin{array}{rcl}
			S_{h,g}f&=&\displaystyle{\sum_{l\in \mathbb{N}_L}\sum_{n\in \mathbb{Z}}\sum_{m\in \mathbb{N}_M}\langle f, E_{\frac{m}{M}}T_{nN}h_l\rangle E_{\frac{m}{M}}T_{nN}g_l}\\
			&=& \displaystyle{\sum_{l\in \mathbb{N}_L}\sum_{n\in \mathbb{Z}}\sum_{m\in \mathbb{N}_M}\sum_{r=1}^R\langle f_r, E_{\frac{m}{M}}T_{nN}h_{l,r}\rangle E_{\frac{m}{M}}T_{nN}g_l}\\
			&=& \displaystyle{\sum_{r=1}^R (  \sum_{l\in \mathbb{N}_L}\sum_{n\in \mathbb{Z}}\sum_{m\in \mathbb{N}_M}\langle f_r, E_{\frac{m}{M}}T_{nN}h_{l,r}\rangle E_{\frac{m}{M}}T_{nN}g_{l,1},\ldots}\\
			&\ldots&, \displaystyle{\sum_{l\in \mathbb{N}_L}\sum_{n\in \mathbb{Z}}\sum_{m\in \mathbb{N}_M}\langle f_r, E_{\frac{m}{M}}T_{nN}h_{l,r} \rangle E_{\frac{m}{M}}T_{nN}g_{l,R})}\\
			&=&\displaystyle{\sum_{r=1}^R\sum_{l\in \mathbb{N}_L}\left( S_{h_{l,r},g_{l,1}}f_r,S_{h_{l,r},g_{l,2}}f_r, \ldots, S_{h_{l,r},g_{l,R}}f_r\right)}.
		\end{array}$$
		Let $k\in \mathbb{N}_p$. We have:
		$$\begin{array}{rcl}
			&&(z_{pM}S_{h,g}f)(j+kM,\theta)\\
			&=&\displaystyle{\sum_{r=1}^R\sum_{l\in \mathbb{N}_L}\left(z_{pM}S_{h_{l,r},g_{l,1}}f_r(j+kM,\theta),\ldots, z_{pM}S_{h_{l,r},g_{l,R}}f_r(j+kM,\theta)\right)}.
		\end{array}$$
		Then,
		$$\begin{array}{rcl}
			&&\left((z_{pM}S_{h,g}f)(j+kM,\theta)\right)_{k\in \mathbb{N}_p}\\
			&=&\displaystyle{\sum_{r=1}^R\sum_{l\in \mathbb{N}_L}\left( \left(z_{pM}S_{h_{l,r},g_{l,1}}f_r(j+kM,\theta)\right)_{k\in \mathbb{N}_p},\ldots, \left( z_{pM}S_{h_{l,r},g_{l,R}}f_r(j+kM,\theta) \right)_{k\in \mathbb{N}_p} \right)}.
		\end{array}$$
		By proposition $2$ in \cite{8}, we have:$$
		\left(\,(z_{pM}S_{h_{l,r},g_{l,r'}}f_r)(j+kM,\theta)\, \right)_{k\in \mathbb{N}_p}
		=M\, Z_{g_{l,r'}}^t(j,\theta) \overline{Z_{h_{l,r}}(j,\theta)} F_r(j,\theta).$$
		Then, Lemmas \ref{lem2.5} and \ref{lem2.6} imply that: 
		$$\begin{array}{rcl}
			\left((z_{pM}S_{h,g}f)(j+kM,\theta)\right)_{k\in \mathbb{N}_p} 
			&=& \displaystyle{\sum_{r=1}^R \sum_{l\in \mathbb{N}_L} \left( Z_{g_{l,1}}^t(j,\theta), Z_{g_{l,2}}^t(j,\theta), \ldots, Z_{g_{l,R}}^t(j,\theta) \right)} \\
			&& \quad \times \left( I_R \otimes \,(\overline{Z_{h_{l,r}}(j,\theta)} F_r(j,\theta)\,) \right) \\
			&=& \displaystyle{\sum_{l\in \mathbb{N}_L} \left( Z_{g_{l,1}}^t(j,\theta), Z_{g_{l,2}}^t(j,\theta), \ldots, Z_{g_{l,R}}^t(j,\theta) \right)} \\
			&& \quad \times \left( I_R \otimes \,(\displaystyle{\sum_{r=1}^R \overline{Z_{h_{l,r}}(j,\theta)} F_r(j,\theta)\,) }\right) \\
			&=& \displaystyle{\sum_{l\in \mathbb{N}_L} \left( Z_{g_{l,1}}^t(j,\theta), Z_{g_{l,2}}^t(j,\theta), \ldots, Z_{g_{l,R}}^t(j,\theta) \right)} \\
			&& \quad \times \left( I_R \otimes \,(\overline{Z_{h_l}(j,\theta)} F(j,\theta)\,) \right).
		\end{array}$$
	\end{proof}
	
	\begin{proposition}\label{prop3.2}
		Let $g:=\{g_l\}_{l\in \mathbb{N}_L},h:=\{h_l\}_{l\in \mathbb{N}_L}\subset \ell^2(\mathbb{Z})$ such that $\mathcal{G}(g,L,M,N,R)$ and $\mathcal{G}(h,L,M,N,R)$ are both Bessel sequences. Then,  the following statements are equivalent.
		\begin{enumerate}
			\item $\mathcal{G}(g,L,M,N,R)$ and $\mathcal{G}(h,L,M,N,R)$ are strongly disjoint (\cite{2}).
			\item $\displaystyle{\sum_{l\in \mathbb{N}_L}Z_{g_l}^*(j,\theta)Z_{h_l}(j,\theta)=0}$ for $j\in \mathbb{N}_{\frac{M}{q}}$ and a.e $\theta\in [0,1[$.\\
		\end{enumerate}
	\end{proposition}
	
	\begin{proof}
		We know that $\mathcal{G}(g,L,M,N,R)$ and $\mathcal{G}(h,L,M,N,R)$ are strongly disjoint if and only if $S_{h,g}f=0$ for all $f\in \ell^2(\mathbb{Z})$ (\cite{2}). Let $f\in \ell^2(\mathbb{Z})$. 
		By choosing $R=1$ in Proposition $\ref{prop1}$, we obtain that $\mathcal{G}(g,L,M,N,R)$ and $\mathcal{G}(h,L,M,N,R)$ are strongly disjoint if and only if $\displaystyle{\sum_{l\in \mathbb{N}_l} Z_{g_l}^t(j,\theta)\overline{Z_{h_l}(j,\theta)}F(j,\theta)=0}$ for $j\in \mathbb{N}_M$ and a.e $\theta \in [0,1[$. Then, $\mathcal{G}(g,L,M,N,R)$ and $\mathcal{G}(h,L,M,N,R)$ are strongly disjoint if and only if $\displaystyle{\sum_{l\in \mathbb{N}_l} Z_{g_l}^*(j,\theta)Z_{h_l}(j,\theta)\overline{F(j,\theta)}=0}$ for $f\in \ell^2(\mathbb{Z})$, $j\in \mathbb{N}_M$ and a.e $\theta \in [0,1[$. For an arbitrary $x\in \mathbb{C}^p$, define $f\in \ell^2(\mathbb{Z})$ by $\overline{F(j,\theta)}=x$ for $j\in \mathbb{N}_M$ and a.e $\theta\in [0,1[$. Then, $\displaystyle{\sum_{l\in \mathbb{N}_l} Z_{g_l}^*(j,\theta)Z_{h_l}(j,\theta)x=0}$ for $j\in \mathbb{N}_M$ and a.e $\theta \in [0,1[$. Hence, \begin{equation}
			\displaystyle{\sum_{l\in \mathbb{N}_l} Z_{g_l}^*(j,\theta)Z_{h_l}(j,\theta)}=0 \text{ for } j\in \mathbb{N}_M \text{ and a.e } \theta\in  [0,1[.
		\end{equation} 
		Evidently, $(3)$ implies that $\displaystyle{\sum_{l\in \mathbb{N}_l} Z_{g_l}^*(j,\theta)Z_{h_l}(j,\theta)\overline{F(j,\theta)}=0}$ for $f\in \ell^2(\mathbb{Z})$, $j\in \mathbb{N}_M$ and a.e $\theta \in [0,1[$, then that $\mathcal{G}(g,L,M,N,R)$ and $\mathcal{G}(h,L,M,N,R)$ are strongly disjoint. Let's prove now that $(3)$ is equivalent to:
		\begin{equation}
			\displaystyle{\sum_{l\in \mathbb{N}_l} Z_{g_l}^*(j,\theta)Z_{h_l}(j,\theta)}=0 \text{ for } j\in \mathbb{N}_{\frac{M}{q}} \text{ and a.e } \theta\in  [0,1[.
		\end{equation} 
		For this,  we prove only that $(4)$ implies $(3)$. By Lemma \ref{lem2.12}, we have for all $l\in \mathbb{N}_L$, $\ell\in \mathbb{N}_q, \; j\in \mathbb{N}_{\frac{M}{q}}$ and a.e $\theta \in [0,1[$,
		$$Z_{g_l}^*(j+\displaystyle{\frac{M}{q}}\ell,\theta)Z_{h_l}(j+\displaystyle{\frac{M}{q}}\ell,\theta)=A_{\ell}^*(\theta)Z_{g_l}^*(j,\theta)Z_{h_l}(j,\theta)A_{\ell}(\theta).$$
		Thus,
		$$\sum_{l\in \mathbb{N}_L} Z_{g_l}^*(j+\displaystyle{\frac{M}{q}}\ell,\theta)Z_{h_l}(j+\displaystyle{\frac{M}{q}}\ell,\theta)=A_{\ell}^*(\theta)\left(\sum_{l\in \mathbb{N}_L} Z_{g_l}^*(j,\theta)Z_{h_l}(j,\theta)\right)A_{\ell}(\theta),$$
		for $\ell\in \mathbb{N}_q, \; j\in \mathbb{N}_{\frac{M}{q}}$ and a.e $\theta \in [0,1[$. Hence, $(4)$ implies $(3)$.
	\end{proof}
	\begin{lemma}\label{lem15}
		Let $f\in \ell^2(\mathbb{Z})$. If $f\in \ell^2(\mathbb{S})$. Then,
		$$Z_f(j,\theta)\mathcal{K}(j)=Z_f(j,\theta) \text{ for } j\in \mathbb{N}_{\frac{M}{q}} \text{ and a.e }\theta\in [0,1[.$$
	\end{lemma}
	\begin{proof}
		Let $j\in \mathbb{N}_{\frac{M}{q}}$, a.e $\theta \in [0,1[$. Let $s\in \mathbb{N}_q$ and $t\in \mathbb{N}_p$. We have,
		$$\begin{array}{rcl}
			\left(Z_f(j,\theta)\mathcal{K}(j)\right)_{s,t}&=&\displaystyle{\sum_{k=0}^{p-1}Z_f(j,\theta)_{s,k}\mathcal{K}(j)_{k,t}}\\
			&=&\displaystyle{\sum_{k=0}^{p-1}Z_f(j,\theta)_{s,k}\delta_{k,t}\chi_{\mathcal{K}_j}(t)}\\
			&=&Z_f(j,\theta)_{s,t}\chi_{\mathcal{K}_j}(t)\\
			&=&\left\lbrace
			\begin{array}{rcl}
				&Z_f(j,\theta)_{s,t}& \; \text{ if } t\in \mathcal{K}_j,\\
				&0& \;\; \text{ otherwise}. 
			\end{array}
			\right.
		\end{array}$$
		On the other hand, we have $Z_f(j,\theta)_{s,t}=z_{pM}f(j+tM-sN,\theta)=\displaystyle{\sum_{k\in \mathbb{Z}}f(j+tM-sN+kpM)e^{2\pi i k\theta)}}=\displaystyle{\sum_{k\in \mathbb{Z}}f(j+tM-sN+kqN)e^{2\pi i k\theta}}$ since $pM=qN$. Then, if $t\notin \mathcal{K}_j$, $j+tM\notin \mathbb{S}$, then, for all $k\in \mathbb{Z}$, $j+tM-sN+kqN\notin \mathbb{S}$ by the $N\mathbb{Z}$-periodicity of $\mathbb{S}$, thus $f(j+tM-sN+kqN)=0$ for all $k\in \mathbb{Z}$. Hence $Z_f(j,\theta)_{s,t}=0$ if $t\notin \mathcal{K}_j$. The proof is completed.
	\end{proof}

	The following proposition shows the link between the rank of $Z_g(j,\theta)$  and cardinality of $\mathcal{K}_j$.
	\begin{proposition}\label{prop3}
		Let $g:=\{g_l\}_{l\in \mathbb{N}_L}\subset \ell^2(\mathbb{S},\mathbb{C}^R)$. Then, for $j\in \mathbb{N}_{\frac{M}{q}}$ and a.e $\theta\in [0,1[$, we have,
		\begin{equation}
			rank(Z_g(j,\theta))\leq R.card(\mathcal{K}_j).
		\end{equation}
	\end{proposition}
	\begin{proof}
		Let $j\in \mathbb{N}_{\frac{M}{q}}$ and a.e $\theta\in [0,1[$. By Lemma \ref{lem15}, we have, for all $l\in \mathbb{N}_L$ and $1\leq r \leq R$, $Z_{g_{l,r}}(j,\theta)\mathcal{K}(j)=Z_{g_{l,r}}(j,\theta)$, where for all $l\in \mathbb{N}_L$, $g_l:=(g_{l,1},g_{l,2},\ldots,g_{l,R})$. Then, by Lemma $\ref{lem2.5}$, we have, for all $j\in \mathbb{N}_{\frac{M}{q}}$ and a.e $\theta\in  [0,1[$,
		\begin{equation}
			Z_g(j,\theta)\left( I_R\otimes \mathcal{K}(j)\right)=Z_g(j,\theta).
		\end{equation} 
		Thus, $$rank(Z_g(j,\theta)\,)=rank\left(Z_g(j,\theta)\left( I_R\otimes \mathcal{K}(j)\right)\, \right)\leq rank( I_R\otimes \mathcal{K}(j))\leq R.card(\mathcal{K}_j).$$
	\end{proof}
	
	\begin{remark}\label{rem2}
		The inequality $(5)$ in Proposition \ref{prop3} holds for all $j\in \mathbb{Z}$ and a.e $\theta\in [0,1[$.\\
		In fact, Lemma \ref{lem2.12} implies that $rank\left(Z_g(j+\displaystyle{\frac{M}{q}}\ell)\right)=rank(Z_g(j,\theta)\,)$ for all $j\in \mathbb{N}_{\frac{M}{q}}$, $\ell\in \mathbb{Z}$ and a.e $\theta\in [0,1[$ since $A_{\ell}(\theta)$ and $C_{\ell}(\theta)$ are both unitaries. And by Remark $2.5$ in \cite{4}, $card(\mathcal{K}_j)$ is $\displaystyle{\frac{M}{q}}$-periodic with respect to $j$. \\
	\end{remark}
	
	The following theorem characterizes which $g$ generates a complete super Gabor system on $\mathbb{S}$.
	\begin{theorem}\label{prop4}
		Let $g:=\{g_l\}_{l\in \mathbb{N}_L}\subset \ell^2(\mathbb{S},\mathbb{C}^R)$. The following statements are equivalent.
		\begin{enumerate}
			\item $\mathcal{G}(g,L,M,N,R)$ is complete in $\ell^2(\mathbb{S},\mathbb{C}^R)$.
			\item For all $j\in \mathbb{N}_{\frac{M}{q}}$ and a.e $\theta \in [0,1[$, \begin{equation}
				rank(Z_g(j,\theta))=R.card(\mathcal{K}_j).
			\end{equation}
		\end{enumerate}
	\end{theorem}
	
	For the proof, we will need the following lemma.
	\begin{lemma}\label{lem14} Let $f\in \ell^2(\mathbb{Z})$. If $f\in \ell^2(\mathbb{S})$.  Then,
		$$\mathcal{K}(j)F(j,\theta)=F(j,\theta)\text{ for }j\in \mathbb{N}_M \text{ and a.e }\theta\in [0,1[.$$
	\end{lemma}
	\begin{proof}
		Let $j\in \mathbb{N}_M$, a.e $\theta \in [0,1[$.  Let $k\in \mathbb{N}_p$. We have, $$\begin{array}{rcl}
			\left( \mathcal{K}(j)F(j,\theta)\right)_k&=&\displaystyle{\sum_{n\in \mathbb{N}_p} (\,\mathcal{K}(j)\,)_{k,n} F(j,\theta)_n}\\
			&=& \displaystyle{\sum_{n\in \mathbb{N}_p} \delta_{n,k} \chi_{\mathcal{K}_j}(k) F(j,\theta)_n}\\
			&=&\chi_{\mathcal{K}_j}(k) F(j,\theta)_k\\
			&=&\left\lbrace
			\begin{array}{rcl}
				&F(j,\theta)_k& \; \text{ if } k\in \mathcal{K}_j,\\
				&0& \;\; \text{ otherwise}. 
			\end{array}
			\right.
		\end{array}$$
		On the other hand, we have $F(j,\theta)_k=z_{pM}f(j+kM,\theta)=\displaystyle{\sum_{n\in \mathbb{Z}}f(j+kM+npM,\theta)e^{2\pi i n\theta}}.$ Then if $k\notin \mathcal{K}_j$ (i.e $j+kM\notin \mathbb{S}$), $j+kM+npM\notin \mathbb{S}$, for all $n\in \mathbb{Z}$, by $N\mathbb{Z}$-periodicity of $\mathbb{S}$ and since $pM=qN$. Hence, if $k\notin \mathcal{K}_j$, $F(j,\theta)_k=0$. Hence, $\left( \mathcal{K}(j)F(j,\theta)\right)_k=F(j,\theta)_k$ for all $k\in \mathbb{N}_p$. Thus, $\mathcal{K}(j)F(j,\theta)=F(j,\theta)$.
	\end{proof}
	\begin{proof}[Proof of Theorem $\ref{prop4}$]
		By Lemma \ref{lem13}, $f\in \ell^2(\mathbb{S},\mathbb{C}^R)$ is orthogonal to $\mathcal{G}(g,L,M,N)$, if and only if, $Z_g(j,\theta)\overline{F(j,\theta)}=0$ for $j\in \mathbb{N}_{M}$ and a.e $\theta\in [0,1[$, which, together with Lemmas \ref{lem15} and \ref{lem14}, is equivalent to:
		\begin{equation}
			(\, Z_g(j,\theta)(I_R\otimes \mathcal{K}(j)\,)\,)(\,(I_R\otimes \mathcal{K}(j))\overline{F(j,\theta)})=0 \text{ for } j\in \mathbb{N}_{M} \text{ and a.e } \theta \in [0,1[.
		\end{equation}
		By Remark $\ref{rem2}$, we only need to prove that $\mathcal{G}(g,L,M,N,R)$ is complete, if and only if, 
		\begin{equation}
			rank( Z_g(j,\theta)(I_R\otimes \mathcal{K}(j))=R.card(\mathcal{K}_j) \text{ for } j\in \mathbb{N}_M \text{ and a.e } \theta \in [0,1[.
		\end{equation}
		Suppose $(9)$ and let $f\in \ell^2(\mathbb{S},\mathbb{C}^R)$ orthogonal to $\mathcal{G}(g,L,M,N,R)$. $(9)$ Shows that the rank of $Z_g(j,\theta)(I_R\otimes \mathcal{K}(j)\,)$ is exactly the number of nonzero columns of $Z_g(j,\theta)(I_R\otimes \mathcal{K}(j)\,)$ for $j\in \mathbb{N}_M$ and a.e $\theta \in [0,1[$. Then, from $(8)$, it follows that $(\,(I_R\otimes \mathcal{K}(j))\overline{F(j,\theta)})=0$, and thus by Lemma \ref{lem15}, $F(j,\theta)=0$ for $j\in \mathbb{N}_M$ and a.e $\theta \in [0,1[$. Hence $f=0$. Hence $\mathcal{G}(g,L,M,N,R)$ is complete in $\ell^2(\mathbb{S},\mathbb{C}^R)$. Conversely, suppose that $\mathcal{G}(g,L,M,N,R)$ is complete in $\ell^2(\mathbb{S},\mathbb{C}^R)$ and suppose that $(9)$ fails. Then there exists $j_0\in \mathbb{N}_M$ and $E_0\subset [0,1[$ with positive measure such that for all $\theta \in E_0$, 
		\begin{equation}
			\begin{array}{rcl}
				rank(Z_g(j_0,\theta))&=&rank(Z_g(j_0,\theta)(I_R\otimes \mathcal{K}(j))\\
				&<& R. card(\mathcal{K}_{j_0})\\
				&=&card(\Lambda_{j_0})
			\end{array}
		\end{equation}
		Denote by $\mathbb{P}(j_0,\theta):\mathbb{C}^{pR}\rightarrow \mathbb{C}^{pR}$ the orthogonal projection onto the kernel of $Z_g(j_0,\theta)$, $N\left( Z_g(j_0,\theta)\right)$, for a.e $\theta\in [0,1[$. Let $\{e_k\}_{k\in \mathbb{N}_{pR}}$ be the standard orthonormal basis for $\mathbb{C}^{pR}$. Suppose that $span\{e_k:\; k\in \Lambda_{j_0}\}\subset N(\mathbb{P}(j_0,\theta))$ for some $\theta \in E_0$. Then $\{e_k:\; k\in \Lambda_{j_0}\}\oplus N(Z_g(j_0,\theta))$ is an orthogonal sum, then: $$\begin{array}{rcl}
			pR&\geqslant& card(\Lambda_{j_0})+(pR-rank(Z_g(j_0,\theta))\,)\\
		\end{array}$$
		Then: $$rank(Z_g(j_0,\theta))\geqslant card(\Lambda_{j_0}).$$
		Contradiction with $(10)$. Therefore, there exist $k_0\in \Lambda_{j_0}$ and $\tilde{E_0}\subset E_0$ with positive measure such that $e_{k_0}\notin N(\mathbb{P}(j_0,\theta))$ for all $\theta \in \tilde{E_0}$. That means that $\mathbb{P}(j_0,\theta)e_{k_0}\neq 0$ for all $\theta \in \tilde{E_0}$. Define $f\in  \ell^2(\mathbb{Z},\mathbb{C}^R)$ such that $F(j,\theta):=\delta_{j,j_0}\mathbb{P}(j_0,\theta)e_{k_0}$ for all $j\in \mathbb{N}_M$ and a.e $\theta\in [0,1[$. Then $\|F(j,\theta)\|_{\mathbb{C}^{pR}}\leq 1$, and the $k$-th component of $F(j_0,\theta)$ equals $\langle \mathbb{P}(j_0,\theta)e_{k_0},e_k\rangle=\langle e_{k_0},\mathbb{P}(j_0,\theta)e_k\rangle$ for $k\in \mathbb{N}_{pR}$ and a.e $\theta\in [0,1[$. Observe that when $k\notin \Lambda_{j_0},\; e_k\in N(Z_g(j_0,\theta))$ for a.e $\theta \in [0,1[$ by the definition of $z_{pM}$, which implies that the $k$-th component of $F(j_0,\theta)$ equals $\langle e_{k_0},\mathbb{P}(j_0,\theta)e_k\rangle=\langle e_{k_0},e_k\rangle=0$ for a.e $\theta\in [0,1[$. Then, by Lemma \ref{lem2.1}, $f\in \ell^2(\mathbb{S},\mathbb{C}^R)$. Since $F(j_0,\theta)=\mathbb{P}(j_0,\theta)e_{k_0}\neq 0$ for all $\theta\in \tilde{E_0}$ which is with positive measure, then $f\neq 0$. On the other hand, we have for $j\neq j_0$, $F(j,\theta)=0$, then $Z_g(j,\theta)F(j,\theta)=0$ for a.e $\theta \in [0,1[$ and $Z_g(j_0,\theta)F(j_0,\theta)=Z_g(j_0,\theta)\mathbb{P}(j_0,\theta)e_{k_0}=0$ for a.e $\theta \in [0,1[$ since $\mathbb{P}(j_0,\theta)e_{k_0}\in N(\mathbb{P}(j_0,\theta))$. Hence, by Lemma $\ref{lem13}$,  $f$ is orthogonal to $\mathcal{G}(g,L,M,N,R)$, but $f\neq 0$. Contradiction with the fact that $\mathcal{G}(g,L,M,N,R)$ is complete in $\ell^2(\mathbb{S},\mathbb{C}^R)$. 
	\end{proof}

	\begin{remark}
		For  the special  case of $\mathbb{S}=\mathbb{Z}$. For all $j\in \mathbb{N}_{\frac{M}{q}}$, $\mathcal{K}_j=\mathbb{N}_p$. Then the condition $(2)$ in the Theorem \ref{prop4} is equivalent to $rank(Z_g(j,\theta)\,)=pR$ for all $j\in \mathbb{N}_{\frac{M}{q}}$ and a.e $\theta \in [0,1[$, which is equivalent to: $Z_g(j,\theta)$ is injective for all $j\in \mathbb{N}_{\frac{M}{q}}$ and a.e $\theta \in [0,1[$.\\
	\end{remark}
	
	The following theorem characterizes which $g$ generates a multi-window super Gabor frame on $\mathbb{S}$.
	\begin{theorem}\label{prop5}
		Given \(g:=\{g_l\}_{l\in \mathbb{N}_L} \subset \ell^2(\mathbb{S},\mathbb{C}^R)\). Then the following statements are equivalent:
		\begin{enumerate}
			\item $\mathcal{G}(g,L, M,N,R)\) is a frame for \(\ell^2(\mathbb{S},\mathbb{C}^R)$ with frame bounds $0 < A \leq B $.
			\item 
			\begin{equation}
				\frac{A}{M}.I_R\otimes \mathcal{K}(j)
				\leq \sum_{l\in \mathbb{N}_L} Z_{g_l}^*(j, \theta)Z_{g_l}(j, \theta) \leq \frac{B}{M}.I_R\otimes \mathcal{K}(j).
			\end{equation}
			for all $j \in \mathbb{N}_{\frac{M}{q}}$ and a.e $\theta \in [0, 1[$.
			\item The inequality $(11)$ holds for all $j\in \mathbb{N}_M$ and a.e $\theta \in [0,1[$.
		\end{enumerate}
	\end{theorem}
	For the proof, we will need the following lemma.
	\begin{lemma}\label{lem16}
		Denote $L^{\infty}(\mathbb{S}_{pM}\times [0,1[,\mathbb{C}^{R})$ the set of vector-valued functions $H$ on $\mathbb{S}_{pM}\times [0,1[$ such that for all $j\in \mathbb{S}_{pM}$, $H(j,.)\in L^{\infty}([0,1[,\mathbb{C}^{R})$, and  $\Delta := z_{pM \mid \ell^2(\mathbb{S},\mathbb{C}^R)}^{-1}\left(L^{\infty}(\mathbb{S}_{pM} \times [0,1[,\mathbb{C}^R)\right).$ Let $g:=\{g_l\}_{l\in \mathbb{N}_L}\subset \ell^2(\mathbb{S},\mathbb{C}^{R})$. Then the following statements are equivalent:
		\begin{enumerate}
			\item $\mathcal{G}(g,L,M,N,R)$ is a frame for $\ell^2(\mathbb{S},\mathbb{C}^{R})$ with frame bounds $A\leq B$.
			\item For all  $f \in \Delta$, we have,
			$$
			\begin{array}{rcl}
				\displaystyle{\frac{A}{M} \sum_{j=0}^{M-1} \int_{0}^{1} \|F(j, \theta)\|^2  d\theta} &\leq& \displaystyle{\sum_{l=0}^{L-1}\sum_{j=0}^{M-1} \int_{0}^{1} \| Z_{g_l}(j, \theta) F(j, \theta)\|^2  d\theta}\\
				& \leq& \displaystyle{\frac{B}{M} \sum_{j=0}^{M-1} \int_{0}^{1} \|F(j, \theta)\|^2 d\theta.}
			\end{array}$$
			
		\end{enumerate}
	\end{lemma}
	\begin{proof}
		By density of $L^{\infty}(\mathbb{S}_{pM}\times [0,1[,\mathbb{C}^R)$ in $L^2(\mathbb{S}_{pM}\times [0,1[,\mathbb{C}^R)$ and by the unitarity of $z_{pM}$ from $\ell^2(\mathbb{S},\mathbb{C}^R)$ onto $L^2(\mathbb{S}_{pM}\times [0,1[,\mathbb{C}^R)$ (Since $pM=qN$, then $\mathbb{S}$ is $pM\mathbb{Z}$-periodic in $\mathbb{Z}$), $\Delta$ is dense in $\ell^2(\mathbb{S},\mathbb{C}^R)$. Hence $\mathcal{G}(g,L,M,N,R)$ is a frame for $\ell^2(\mathbb{S},\mathbb{C}^R)$ with frame bounds $A\leq B$, if and only if, for all $f\in \Delta$, $$A\|f\|^2\leq \displaystyle{\sum_{l\in \mathbb{N}_L}\sum_{n\in \mathbb{Z}}\sum_{m\in \mathbb{N}_M}\left\vert \langle f, E_{\frac{m}{M}}T_{nN}g_l\rangle\right\vert^2}\leq B\|f\|^2,$$ if and only if, for all $f\in \Delta$, $$A\|f\|^2\leq \displaystyle{\sum_{l\in \mathbb{N}_L }\sum_{r\in \mathbb{N}_q}\sum_{n\in \mathbb{Z}}\sum_{m\in \mathbb{N}_M}\left\vert \langle f, E_{\frac{m}{M}}T_{(nq+r)N}g_l\rangle\right\vert^2}\leq B\|f\|^2.$$
		Let $f\in \Delta$, by Lemma \ref{lem2.4} we have, 
		$$
		\begin{array}{rcl}
			\langle f, E_{\frac{m}{M}}T_{(nq+r)N}g_l\rangle &=&\displaystyle{\sum_{j\in \mathbb{N}_M}\left( \int_0^1 (\overline{Z_{g_l}(j,\theta)}F(j,\theta))_re^{-2\pi i n\theta}d\theta\right) e^{-2\pi i \frac{m}{M}j}}\\
			&=& \displaystyle{\sum_{j\in \mathbb{N}_M}T(j)e^{-2\pi i \frac{m}{M}j}}\\
			&=&\langle T,e^{2\pi i \frac{m}{M}.}\rangle,
		\end{array}$$
		where  $T(j)=\int_0^1 (\overline{Z_{g_l}(j,\theta)}F(j,\theta))_re^{-2\pi i n\theta}d\theta$. Observe that $T$ is $M$-periodic. Since $\{\displaystyle{\frac{1}{\sqrt{M}}e^{2\pi i \frac{m}{M}.}}\}_{m\in \mathbb{N}_M}$ is an orthonormal basis for $\ell^2(\mathbb{N}_M)$; the space of $M$-periodic sequences, then we have:$$\begin{array}{rcl}
			&&\displaystyle{\sum_{m\in \mathbb{N}_M}\vert \langle f, E_{\frac{m}{M}}T_{(nq+r)N}g_l\rangle\vert^2}\\
			&=& \displaystyle{\sum_{m\in \mathbb{N}_M}\left \vert \langle T, e^{2\pi i \frac{m}{M}.}\rangle \right\vert^2}\\
			&=& M\|T\|^2\\
			&=&\displaystyle{M\sum_{j\in \mathbb{N}_M}\left\vert\int_0^1 (\overline{Z_{g_l}(j,\theta)}F(j,\theta)\,)_{r}\,e^{-2\pi i n\theta}d\theta \right\vert^2}.\\
		\end{array}$$
		Since $\{e^{2\pi i n\theta}\}_{n\in \mathbb{Z}}$ is an orthonormal basis for $L^2([0,1[)$ and $ (\overline{Z_{g_l}(j,\theta)}F(j,\theta)\,)_{r}\in L^2([0,1[)$ for all $r\in \mathbb{N}_q$ , then we have: 
		$$\begin{array}{rcl}
			&&\displaystyle{\sum_{m\in \mathbb{N}_M}\sum_{n\in \mathbb{Z}}\vert \langle f, E_{\frac{m}{M}}T_{(nq+r)N}g_l\rangle\vert^2}\\
			&=&\displaystyle{M\sum_{j\in \mathbb{N}_M}\sum_{n\in \mathbb{Z}}\left\vert\int_0^1 (\overline{Z_{g_l}(j,\theta})F(j,\theta))_{r}\,e^{-2\pi i n\theta}d\theta \right\vert^2}\\
			&=&\displaystyle{M\sum_{j\in \mathbb{N}_M}\sum_{n\in \mathbb{Z}}\left\vert \langle (\overline{Z_{g_l}(j,.)}F(j,.))_{r},e^{2\pi i n.}\rangle \right\vert^2}\\
			&=&\displaystyle{M\sum_{j\in \mathbb{N}_M}\left\| (\overline{Z_{g_l}(j,.)}F(j,.))_r\right\|^2}\\
			&=&\displaystyle{M\sum_{j\in \mathbb{N}_M}\int_0^1 \left\vert (Z_{g_l}(j,\theta)F(j,\theta))_r\right\vert^2\,d\theta}\\
		\end{array}$$
		Hence, $$\begin{array}{rcl}
			&&\displaystyle{\sum_{r\in \mathbb{N}_q}\sum_{n\in \mathbb{Z}}\sum_{m\in \mathbb{N}_M}\vert \langle f, E_{\frac{m}{M}}T_{(nq+r)N}g_l\rangle\vert^2}\\
			&=&\displaystyle{M\sum_{j\in \mathbb{N}_M}\int_0^1 \sum_{r\in \mathbb{N}_q}\left\vert (Z_{g_l}(j,\theta)F(j,\theta))_r\right\vert^2\,d\theta}\\
			&=&\displaystyle{M\sum_{j\in \mathbb{N}_M}\int_0^1 \left\| Z_{g_l}(j,\theta)F(j,\theta)\right\|^2\,d\theta}.\\
		\end{array}$$
		The norm in the last line is the 2-norm in $\mathbb{C}^q$. Thus,
		$$\begin{array}{rcl}
			&&\displaystyle{\sum_{l\in \mathbb{N}_L}\sum_{r\in \mathbb{N}_q}\sum_{n\in \mathbb{Z}}\sum_{m\in \mathbb{N}_M}\vert \langle f, E_{\frac{m}{M}}T_{(nq+r)N}g_l\rangle\vert^2}\\ &=&\displaystyle{M\sum_{l\in \mathbb{N}_L}\sum_{j\in \mathbb{N}_M}\int_0^1 \left\| Z_{g_l}(j,\theta)F(j,\theta)\right\|^2\,d\theta}.
		\end{array}$$
		Hence, \begin{equation}
			\displaystyle{\sum_{l\in \mathbb{N}_L}\sum_{r\in \mathbb{N}_q}\sum_{n\in \mathbb{Z}}\sum_{m\in \mathbb{N}_M}\vert \langle f, E_{\frac{m}{M}}T_{(nq+r)N}g_l\rangle\vert^2}=\displaystyle{M\sum_{l\in \mathbb{N}_L}\sum_{j\in \mathbb{N}_M}\int_0^1 \left\| Z_{g_l}(j,\theta)F(j,\theta)\right\|^2\,d\theta}.
		\end{equation}
		On the other hand, we have by unitarity of $z_{pM}$:$$\begin{array}{rcl}
			\|f\|^2&=&\|z_{pM}f\|^2\\
			&=&\displaystyle{\sum_{r\in \mathbb{N}_R}\|z_{pM}f_r\|^2}\\
			&=&\displaystyle{\sum_{r\in \mathbb{N}_R}\sum_{j\in \mathbb{N}_{pM}}\int_0^1\vert z_{pM}f_r(j,\theta)\vert^2\, d\theta}\\
			&=&\displaystyle{\sum_{r\in \mathbb{N}_R}\sum_{j\in \mathbb{N}_M}\sum_{k\in \mathbb{N}_p}\int_0^1 \vert z_{pM}f_r(j+kM,\theta)\vert^2\, d\theta}\\
			&=&\displaystyle{\sum_{j\in \mathbb{N}_M}\int_0^1 \sum_{r\in \mathbb{N}_R}\sum_{k\in \mathbb{N}_p} \vert F_r(j,\theta)_k\vert^2\, d\theta}\\
			&=&\displaystyle{\sum_{j\in \mathbb{N}_M}\int_0^1 \|F(j,\theta)\|^2\, d\theta}\\
		\end{array}$$
		The norm in the last line is the 2-norm in $\mathbb{C}^{pR}$. Thus, \begin{equation}
			\|f\|^2=\displaystyle{\sum_{j\in \mathbb{N}_M}\int_0^1 \| F(j,\theta)\|^2\, d\theta}.
		\end{equation}
		Then, combining $(12)$ and $(13)$, the proof is completed.
	\end{proof}
	\begin{proof}[Proof of Theorem \ref{prop5}]\hspace{1cm}
		\begin{enumerate}
			\item[]$1.\Longrightarrow 3. $: Assume that $\mathcal{G}(g,L,M,N,R)$ is a frame for $\ell^2(\mathbb{S})$. Then for all $f\in \Delta$, \begin{equation*}
				\begin{array}{rcl}
					\displaystyle{\frac{A}{M} \sum_{j=0}^{M-1} \int_{0}^{1} \|F(j, \theta)\|^2 \, d\theta} &\leq& \displaystyle{\sum_{l=0}^{l-1}\sum_{j=0}^{M-1} \int_{0}^{1} \| Z_{g_l}(j, \theta) F(j, \theta)\|^2 \, d\theta } \\
					&\leq& \displaystyle{\frac{B}{M} \sum_{j=0}^{M-1} \int_{0}^{1} \|F(j, \theta)\|^2 \, d\theta.}
				\end{array}	
			\end{equation*}
			Fix $x:=\{x_k\}_{k\in \mathbb{N}_{pR}}\in  \mathbb{C}^{pR}$, $j_0\in \mathbb{N}_M$ and $h\in L^{\infty}([0,1[)$, and define for all $j\in \mathbb{N}_M$ and a.e $\theta\in [0,1[$, $F(j,\theta):=\{\delta_{j,j_0}\,\chi_{\Lambda_j}(k)\,x_k\, h(\theta) \}_{k\in \mathbb{N}_{pR}}$. \\Then, $\displaystyle{ \sum_{j=0}^{M-1} \int_{0}^{1} \|F(j, \theta)\|^2 \, d\theta=\|\{\chi_{\Lambda_{j_0}}(k)x_k\}_{k\in \mathbb{N}_{pR}}\|^2\int_0^1 \vert h(\theta)\vert^2 \,d\theta}$. Since 
			$\{\chi_{\Lambda_{j_0}}(k)x_k\}_{k\in \mathbb{N}_{pR}}=diag\left(\chi_{\Lambda_{j_0}}(0),\chi_{\Lambda_{j_0}}(1),\ldots,\chi_{\Lambda_{j_0}}(pR-1)\right)x$, then by Lemma \ref{lem2.15}, $\{\chi_{\Lambda_{j_0}}(k)x_k\}_{k\in \mathbb{N}_{pR}}=(I_R\otimes \mathcal{K}(j_0)\,)x$. Thus, by Lemmas \ref{lem2.9}, \ref{lem2.7} and \ref{lem2.14}, we have:
			$\displaystyle{ \sum_{j=0}^{M-1} \int_{0}^{1} \|F(j, \theta)\|^2 \, d\theta}=\displaystyle{\langle (I_R\otimes\mathcal{K}(j_0))x,x\rangle \int_0^1 \vert h(\theta)\vert^2\,d\theta}$. On the other hand, we have: $$
			\begin{array}{rcl}
				&&\displaystyle{\sum_{l=0}^{L-1}\sum_{j=0}^{M-1} \int_{0}^{1} \| Z_{g_l}(j, \theta) F(j, \theta)\|^2 \, d\theta }\\&=&\displaystyle{\sum_{l=0}^{L-1} \int_{0}^{1} \| Z_{g_l}(j_0, \theta) F(j_0, \theta)\|^2 \, d\theta} \\
				&=&\displaystyle{\sum_{l=0}^{L-1} \int_{0}^{1} \sum_{r=0}^{q-1} \left\vert \left(Z_{g_l}(j_0, \theta) F(j_0, \theta)\right)_r \,\right\vert^2 \, d\theta} \\
				&=&\displaystyle{\sum_{l=0}^{L-1} \int_{0}^{1} \sum_{r=0}^{q-1} \left\vert \sum_{k=0}^{pR-1} Z_{g_l}(j_0, \theta)_{r,k} F(j_0, \theta)_k \,\right\vert^2 \, d\theta} \\
				&=&\displaystyle{\sum_{l=0}^{L-1} \int_{0}^{1} \sum_{r=0}^{q-1} \left\vert \sum_{k=0}^{pR-1} Z_{g_l}(j_0, \theta)_{r,k} \,\chi_{\_{j_0}}(k)\,x_k\, h(\theta)\,\right\vert^2 d\theta }\\
				&=&\displaystyle{\sum_{l=0}^{L-1} \int_{0}^{1} \sum_{r=0}^{q-1}\left \vert \sum_{k=0}^{pR-1} Z_{g_l}(j_0, \theta)_{r,k} ((I_R\otimes \mathcal{K}(j_0))x)_k \right\vert^2 \vert h(\theta)\vert^2 d\theta} \\
				&=&\displaystyle{\sum_{l=0}^{L-1} \int_{0}^{1} \sum_{r=0}^{q-1} \left\vert \left(Z_{g_l}(j_0, \theta)(I_R\otimes \mathcal{K}(j_0))x\right)_r \right\vert^2 \vert h(\theta)\vert^2 d\theta}\\
				&=&\displaystyle{\sum_{l=0}^{L-1} \int_{0}^{1} \left\|Z_{g_l}(j_0,\theta)(I_R\otimes\mathcal{K}(j_0))x \right\|^2 \vert h(\theta)\vert^2 d\theta} \\
				&=&\displaystyle{\sum_{l=0}^{L-1} \int_{0}^{1} \left\|Z_{g_l}(j_0,\theta)x \,\right\|^2 \vert h(\theta)\vert^2  d\theta}\;\;\; \text{ by Lemma \ref{lem15}} \\
				&=&\displaystyle{\int_{0}^{1} \left\langle \sum_{l=0}^{L-1} Z_{g_l}^*(j_0,\theta)Z_{g_l}(j_0,\theta)x,x\right\rangle \vert h(\theta)\vert^2 d\theta}. \\
			\end{array}$$
			Then for all $j\in \mathbb{N}_M$,  $x\in \mathbb{C}^{pR}$ and $h\in L^{\infty}([0,1[)$, we have, 
			\begin{equation}
				\begin{array}{rcl}
					\displaystyle{\frac{A}{M}.\left\langle (I_R\otimes\mathcal{K}(j))x,x\right\rangle \,\int_0^1\vert h(\theta)\vert^2 \,d\theta} &\leq& \displaystyle{\int_{0}^{1} \left\langle \sum_{l=0}^{l-1} Z_{g_l}^*(j,\theta)Z_{g_l}(j,\theta)x,x\right\rangle \,\vert h(\theta)\vert^2 \, d\theta} \\
					&\leq& \displaystyle{\frac{B}{M}.\left\langle (I_R\otimes\mathcal{K}(j))x,x\right \rangle \int_0^1\vert h(\theta)\vert^2\, d\theta.}
				\end{array}
			\end{equation}
			For $j\in \mathbb{N}_M$ and $x\in \mathbb{C}^{pR}$ fixed, denote $C=\displaystyle{\frac{A}{M}.\langle (I_R\otimes\mathcal{K}(j))x,x\rangle}$ and $D=\displaystyle{\frac{B}{M}.\langle (I_R\otimes \mathcal{K}(j))x,x\rangle}$. Assume, by contradiction, that \begin{equation} C>\displaystyle{\left\langle \sum_{l=0}^{L-1}Z_{g_l}^*(j,\theta)Z_{g_l}(j,\theta)x,x\right\rangle} ,
			\end{equation} on a subset of $[0,1[$ with a positive measure. Denote 
			$D=\{\theta \in [0,1[:\; (15) \text{ holds}\}.$\\
			For all $k\in \mathbb{N}$, denote: $$D_k:=\left\{\theta\in [0,1[:\; C-\displaystyle{\frac{C}{k+1}}\leq \displaystyle{\langle \sum_{l=0}^{L-1}Z_{g_l}^*(j,\theta)Z_{g_l}(j,\theta)x,x\rangle} \leq C-\displaystyle{\frac{C}{k}}\right\}.$$
			It is clear that $\{D_k\}_{k\in \mathbb{N}}$ forms a partition for $D$. Since $mes(D)>0$, then there exists $k\in \mathbb{N}$ such that $mes(D_k)>0$. Let $h:=\chi_{D_k}$, we have: $$\begin{array}{rcl}
				\displaystyle{\int_{0}^{1} \left\langle \sum_{l=0}^{l-1} Z_{g_l}^*(j,\theta)Z_{g_l}(j,\theta)x,x\right\rangle \,\vert h(\theta)\vert^2 \, d\theta}&=&\displaystyle{\int_{D_k} \left\langle \sum_{l=0}^{l-1} Z_{g_l}^*(j,\theta)Z_{g_l}(j,\theta)x,x\right\rangle d\theta}\\
				&\leq&(C-\displaystyle{\frac{C}{k}}). mes(D_k)\\
				&<&C.mes(D_k)\\
				&=&C\displaystyle{\int_0^1 \vert h(\theta)\vert^2\,d\theta}.
			\end{array}$$	 Contradiction with $(14)$. 
			Suppose, again by contradiction, that \begin{equation} D<\displaystyle{\left\langle \sum_{l=0}^{L-1}Z_{g_l}^*(j,\theta)Z_{g_l}(j,\theta)x,x\right\rangle}, 
			\end{equation} on a subset of $[0,1[$ with a positive measure. Denote 
			$D'=\{\theta \in [0,1[:\; (16) \text{ holds}\}.$\\
			For all $k\in \mathbb{N}$, $m\in \mathbb{N}$, denote $$D_{k,m}':=\left\{\theta\in [0,1[:\; D(k+\displaystyle{\frac{1}{m+1}})\leq \displaystyle{\langle \sum_{l=0}^{L-1}Z_{g_l}^*(j,\theta)Z_{g_l}(j,\theta)x,x\rangle} \leq D(k+\displaystyle{\frac{1}{m}})\right\}.$$
			It is clear that $\{D_{k,m}'\}_{k\in \mathbb{N}}$ forms a partition for $D'$. Since $mes(D')>0$, then there exist $k\in \mathbb{N}$ and $m\in \mathbb{N}$ such that $mes(D_{k,m}')>0$. Let $h:=\chi_{D_{k,m}'}$, we have: $$\begin{array}{rcl}
				\displaystyle{\int_{0}^{1} \langle \sum_{l=0}^{l-1} Z_{g_l}^*(j,\theta)Z_{g_l}(j,\theta)x,x\rangle \,\vert h(\theta)\vert^2 \, d\theta}&=&\displaystyle{\int_{D_{k,m}'} \langle \sum_{l=0}^{l-1} Z_{g_l}^*(j,\theta)Z_{g_l}(j,\theta)x,x\rangle \, d\theta}\\
				&\geqslant&D(k+\displaystyle{\frac{1}{m+1}}). mes(D_{k,m}')\\
				&>&D.mes(D_{k,m}')\\
				&=&D\displaystyle{\int_0^1 \vert h(\theta)\vert^2\,d\theta}.
			\end{array}$$
			Contradiction with $(14)$. Hence, for all $j\in \mathbb{N}_M$, and a.e $\theta \in [0,1[$, we have, 
			\begin{equation*}
				\frac{A}{M}.I_R\otimes\mathcal{K}(j)
				\leq \sum_{l\in \mathbb{N}_L} Z_{g_l}^*(j, \theta)Z_{g_l}(j, \theta) \leq \frac{B}{M}.I_R\otimes\mathcal{K}(j).
			\end{equation*}
			\item[]$(3)\Longrightarrow(1)$: Assume $(3)$. Let $f\in \Delta$, then by Lemma \ref{lem14}, we have for all $j\in \mathbb{N}_M$ and a.e $\theta\in [0,1[$, $(I_R\otimes \mathcal{K}(j)\,)F(j,\theta)=F(j,\theta)$. Then, by $(3)$, we have: \begin{equation*}
				\frac{A}{M}\|F(j,\theta)\|^2\leq \left \langle \sum_{l=0}^{L-1} Z_{g_l}^*(j,\theta)Z_{g_l}(j,\theta)F(j,\theta),F(j,\theta)\right \rangle\leq \frac{B}{M}\|F(j,\theta)\|^2.
			\end{equation*}
			Hence, \begin{equation*}
				\frac{A}{M}\sum_{j=0}^{M-1}\int_0^1\|F(j,\theta)\|^2\leq \sum_{l=0}^{L-1}\sum_{j=0}^{M-1}\int_0^1 \| Z_{g_l}(j,\theta)F(j,\theta)\|^2\leq \frac{B}{M}\sum_{j=0}^{M-1}\int_0^1\|F(j,\theta)\|^2.
			\end{equation*}
			Then, Lemma \ref{lem16} implies $1.$.

			\item[] $2.\Longleftrightarrow 3. $: It is a direct consequence of Lemmas \ref{lem2.12} and \ref{lem2.16}.
		\end{enumerate}
	\end{proof}
	\begin{remark}
		In the special case of $\mathbb{S}=\mathbb{Z}$, for all $j\in \mathbb{N}_{\frac{M}{q}}$, $\mathcal{K}(j)=\mathbb{N}_p$. Then  the condition $(2)$ in the Theorem \ref{prop5} is equivalent to: For all $j\in \mathbb{N}_{\frac{M}{q}}$ and a.e $\theta \in [0,1[$, 
		\begin{equation*}
			\frac{A}{M}.I_{pR}
			\leq \sum_{l\in \mathbb{N}_L} Z_{g_l}^*(j, \theta)Z_{g_l}(j, \theta) \leq \frac{B}{M}.I_{pR}.
		\end{equation*}
		
		Where $I_{pR}$ is the identity matrix in $\mathcal{M}_{pR,pP}$.\\
	\end{remark}
	
	Let $g:=\{g_l\}_{l\in \mathbb{N}_L}\subset \ell^2(\mathbb{S},\mathbb{C}^R)$ and write for all $l\in \mathbb{N}_L$, $g_l=(g_{l,1},g_{l,2},\ldots, g_{l,R})$. Since for all $l\in \mathbb{N}_L$, $j\in \mathbb{N}_{\frac{M}{q}}$ and a.e $\theta\in [0,1[$, we have: $$\begin{array}{rcl}
		Z_{g_l}^*(j,\theta)Z_{g_l}(j,\theta)&=&\begin{pmatrix}
			Z_{g_{l,1}}^*(j,\theta)\\
			Z_{g_{l,2}}^*(j,\theta)\\
			\vdots\\
			Z_{g_{l,R}}^*(j,\theta)
		\end{pmatrix}\times \begin{pmatrix}
			Z_{g_{l,1}}(j,\theta) & Z_{g_{l,2}}(j,\theta)& \ldots&Z_{g_{l,R}}(j,\theta)
		\end{pmatrix}\\
		&=& \left( Z_{g_{l,r}}^*(j,\theta)Z_{g_{l,r'}}(j,\theta)\right)_{1\leq r,r'\leq R},
		\end {array}$$
		where $\left( Z_{g_{l,r}}^*(j,\theta)Z_{g_{l,r'}}(j,\theta)\right)_{1\leq r,r'\leq R}$ is  a block-structured matrix, where $r$ refers to the row index and $r'$ refers to the column index. Hence  for all $j\in \mathbb{N}_{\frac{M}{q}}$ and a.e $\theta \in [0,1[$, we have: 
		$$\displaystyle{\sum_{l\in \mathbb{N}_L}Z_{g_l}^*(j,\theta)Z_{g_l}(j,\theta)=\left( \sum_{l\in \mathbb{N}_L} Z_{g_{l,r}}^*(j,\theta)Z_{g_{l,r'}}(j,\theta)\right)_{1\leq r,r'\leq R}},$$
		where the matrix in the right hand is a block-structured matrix, where $r$ refers to the row index and $r'$ refers to the column index. Thus, Proposition \ref{prop3.2} together with Theorem \ref{prop5} leads to the following two corollaries, as noted as well in \cite{2}.
		\begin{corollary}\label{cor1}
			Let $g:=\{g_l\}_{l\in \mathbb{N}_L}\in \ell^2(\mathbb{S},\mathbb{C}^R)$. Then $\mathcal{G}(g,L,M,N,R)$ is a super Gabor frame for $\ell^2(\mathbb{S},\mathbb{C}^R)$ if  for each $1\leq r\leq R$, $\mathcal{G}(\, \{g_{l,r}\}_{l\in \mathbb{N}_L}, L,M,N)$ is a Gabor frame for $\ell^2(\mathbb{S})$ and $\mathcal{G}(\, \{g_{l,r}\}_{l\in \mathbb{N}_L}, L,M,N)$, for $1\geqslant r\leq R$, are mutually strongly disjoint.\\
		\end{corollary}
		
		\begin{corollary}\label{cor2}
			Let $g:=\{g_l\}_{l\in \mathbb{N}_L}\in \ell^2(\mathbb{S},\mathbb{C}^R)$. Then the following statements are equivalent. 
			\begin{enumerate}
				\item  $\mathcal{G}(g, L,M,N,R)$ is a tight Gabor  frame for $\ell^2(\mathbb{S},\mathbb{C}^R)$ with frame bound $A$.
				\item For each $1\leq r\leq R$, $\mathcal{G}(\, \{g_{l,r}\}_{l\in \mathbb{N}_L}, L,M,N)$ is a Gabor tight frame for $\ell^2(\mathbb{S})$ with frame bound $A$ and $\mathcal{G}(\, \{g_{l,r}\}_{l\in \mathbb{N}_L}, L,M,N)$, for $1\leq r\leq R$, are mutually strongly disjoint.
			\end{enumerate}
		\end{corollary}
		\section{Admissibility Conditions for Multi-window Gabor Systems in \(\ell^2(\mathbb{S}, \mathbb{C}^R)\)
		}
		In this section, we use all the notations already introduced in what above. Let $L,M,N,R\in \mathbb{N}$ and $p,q\in \mathbb{N}$ such that $p\wedge q=1$ and  that $\displaystyle{\frac{N}{M}=\frac{p}{q}}$, and let $\mathbb{S}$ be an $N\mathbb{Z}$-periodic set in $\mathbb{Z}$. In what follows, we give an admissibility condition for $\mathbb{S}$ to admit a complete  Gabor system, a  Gabor (Parseval) frame and a  Gabor (orthonormal) basis. We will need these lemmas.
		\begin{lemma}\label{lem4.2}
			Let $L,M\in \mathbb{N}$ and $E_0,\,E_1,\, \ldots,\,E_{L-1}\subset \mathbb{Z}$ be mutually disjoint. Denote $E=\displaystyle{\bigcup_{l\in \mathbb{N}_L}E_l}$. Then, the following statements are equivalent:
			\begin{enumerate}
				\item[]$(1)$ $\{e^{2\pi i \frac{m}{M}.}\chi_{E_l}\}_{m\in \mathbb{N}_M, \, l\in \mathbb{N}_L}$ is complete in  $\ell^2(E)$.
				\item[] $(2)$ For all $l\in \mathbb{N}_L$, $\{e^{2\pi i \frac{m}{M}.}\chi_{E_l}\}_{m\in \mathbb{N}_M}$ is complete in $\ell^2(E_l)$.
			\end{enumerate}
		\end{lemma}
		\begin{proof}\hspace{1cm}
			\begin{enumerate}
				\item[]$(1)\Longrightarrow (2)$: Assume $(1)$. Fix $l_0\in \mathbb{N}_L$ and let $f\in \ell^2(E_{l_0})$ be orthogonal to  $\{e^{2\pi i \frac{m}{M}.}\chi_{E_{l_0}}\}_{m\in \mathbb{N}_M}$. 
				Define $\overline{f}\in \ell^2(E)$ by $\overline{f}(j)=f(j)$ if $j\in E_{l_0}$ and $0$ otherwise. It is clear that if $l\neq l_0$, $\overline{f}$ is orthogonal to  $\{e^{2\pi i \frac{m}{M}.}\chi_{E_l}\}_{m\in \mathbb{N}_M}$. And we have: 
				$$\begin{array}{rcl}
					\langle \overline{f},e^{2\pi i \frac{m}{M}.}\chi_{E_{l_0}}\rangle&=&\displaystyle{\sum_{j\in E}\overline{f}(j)e^{-2\pi i \frac{m}{M}j}\chi_{E_{l_0}}(j)}\\
					&=& \displaystyle{\sum_{j\in E_{l_0}}f(j)e^{-2\pi i \frac{m}{M}j}\chi_{E_{l_0}}(j)}\\
					&=&0\;\;\; \text{ since } f\text{ is orthogonal to } \{e^{2\pi i \frac{m}{M}.}\chi_{E_{l_0}}\}_{m\in \mathbb{N}_M}.
				\end{array}$$
				Then $\overline{f}$ is orthogonal to $\{e^{2\pi i \frac{m}{M}.}\chi_{E_l}\}_{m\in \mathbb{N}_M, l\in \mathbb{N}_L}$ which is complete in $\ell^2(E)$, thus $\overline{f}=0$ on $E$, and then $f=0$ on $E_l$.
				
				\item[]$(2)\Longrightarrow(1)$: Assume $(2)$ and let $h\in \ell^2(E)$ be orthogonal to $\{e^{2\pi i \frac{m}{M}.}\chi_{E_l}\}_{m\in \mathbb{N}_M, l\in \mathbb{N}_L}$. For all $l\in \mathbb{N}_L$, define $h_l\in \ell^2(E_l)$ as the restriction of $h$ on $E_l$, i.e.  $h_l:=h|_{E_l}$. Let $l\in \mathbb{N}_L$. Fix $l\in \mathbb{N}_L$. Since $h$ is orthogonal to $\{e^{2\pi i \frac{m}{M}.}\chi_{E_l}\}_{m\in \mathbb{N}_M}$, then $\displaystyle{\sum_{j\in E}h(j) e^{-2\pi i \frac{m}{M}.}\chi_{E_l}(j)=0}$, then $\displaystyle{\sum_{j\in E_l}h(j) e^{-2\pi i \frac{m}{M}.}\chi_{E_l}(j)=0}$, thus $\displaystyle{\sum_{j\in E_l}h_l(j) e^{-2\pi i \frac{m}{M}.}\chi_{E_l}(j)=\langle h_l,e^{2\pi i \frac{m}{M}.}\chi_{E_l}\rangle =0}$. Hence $h_l$ is orthogonal to  $\{e^{2\pi i \frac{m}{M}.}\chi_{E_l}\}_{m\in \mathbb{N}_M}$ which is complete in $\ell^2(E_l)$. Hence $h_l=0$ on $E_l$. This for all $l\in \mathbb{N}_L$, therfore $h=0$ on $E$. Hence $\{e^{2\pi i \frac{m}{M}.}\chi_{E_l}\}_{m\in \mathbb{N}_M,l\in \mathbb{N}_L}$ is complete in $\ell^2(E)$.
				
			\end{enumerate}
			
		\end{proof}

		\begin{lemma}\label{lem4.3}
			Let $L,M\in \mathbb{N}$ and $E_0,\,E_1,\, \ldots,\,E_{L-1}\subset \mathbb{Z}$  be mutually disjoint. Denote $E=\displaystyle{\bigcup_{l\in \mathbb{N}_L}E_l}$. Then,  the following statements are equivalent:
			\begin{enumerate}
				\item[] $(1)$ $\{e^{2\pi i \frac{m}{M}.}\chi_{E_l}\}_{m\in \mathbb{N}_M, \, l\in \mathbb{N}_L}$ is a $M$-tight frame for $\ell^2(E)$.
				\item[] $(2)$ For all $l\in \mathbb{N}_L$, $\{e^{2\pi i \frac{m}{M}.}\chi_{E_l}\}_{m\in \mathbb{N}_M}$ is a $M$-tight frame for $\ell^2(E_l)$. 
			\end{enumerate}
		\end{lemma}
		\begin{proof}\hspace{1cm}
			\begin{enumerate}
				\item[]$(1)\Longrightarrow (2)$: Assume $(1)$. Fix $l_0\in \mathbb{N}_L$ and let $f\in \ell^2(E_{l_0})$. Define $\overline{f}\in \ell^2(E)$ by $\overline{f}(j)=f(j)$ if $j\in E_{l_0}$ and $0$ otherwise. It is clear that $\langle \overline{f},e^{2\pi i \frac{m}{M}.}\chi_{E_l}\rangle=0$ if $l\neq l_0$ and that $\|\overline{f}\|=\|f\|$. Together with the fact that $\{e^{2\pi i \frac{m}{M}.}\chi_{E_l}\}_{m\in \mathbb{N}_M, \, l\in \mathbb{N}_L}$ is a tight frame for $\ell^2(E)$ with frame bound $M$, we have: $$
				\begin{array}{rcl}
					M\|f\|^2&=&\displaystyle{\sum_{m=0}^{M-1}\left\vert \sum_{j\in E}\overline{f}(j)e^{-2\pi i \frac{m}{M}j}\chi_{E_{l_0}}(j)\right\vert^2}\\
					&=&\displaystyle{\sum_{m=0}^{M-1}\left\vert\sum_{j\in E_{l_0}}f(j)e^{-2\pi i \frac{m}{M}j}\chi_{E_{l_0}}(j)\right\vert^2}\\
					&=&\displaystyle{\sum_{m=0}^{M-1}\left\vert\langle f,e^{2\pi i \frac{m}{M}.}\chi_{E_{l_0}}\rangle\right\vert^2}.
				\end{array}$$
				Hence, $\{e^{2\pi i \frac{m}{M}.}\chi_{E_{l_0}}\}_{m\in \mathbb{N}_M}$ is a tight frame for $\ell^2(E_{l_0})$ with frame bound $M$. And this for all $l_0\in \mathbb{N}_L$.
				\item[]$(2)\Longrightarrow (1)$: Assume $(2)$. Let $h\in \ell^2(E)$. For all $l\in \mathbb{N}_L$, define $h_l\in \ell^2(E_l)$ as the restriction of $h$ en $E_l$, i.e $h_l=h|_{E_l}$. We have $\displaystyle{\sum_{j\in E}\vert h(j)\vert^2=\sum_{l\in \mathbb{N}_L}\sum_{j\in E_l} \vert h_l(j)\vert^2}$, then $\|h\|^2=\displaystyle{\sum_{l\in \mathbb{N}_L}\|h_l\|^2}.$ It is also clear that  $\langle h,e^{2\pi i \frac{m}{M}.}\chi_{E_l}\rangle=\langle h_l,e^{2\pi i \frac{m}{M}.}\chi_{E_l}\rangle$. Since For all $l\in \mathbb{N}_L$, $\{e^{2\pi i \frac{m}{M}.}\chi_{E_l}\}_{m\in \mathbb{N}_M}$ is a tight frame for $\ell^2(E_l)$ with frame bound $M$, then for all $l\in \mathbb{N}_l$, we have $M\|h_l\|^2=\displaystyle{\sum_{m\in \mathbb{N}_M}\vert \langle h_l,e^{2\pi i \frac{m}{M}.}\chi_{E_l}\rangle \vert^2}.$ Hence, 
				$$M\|h\|^2=\displaystyle{\sum_{l\in \mathbb{N}_L}\sum_{m\in \mathbb{N}_M}\vert \langle h,e^{2\pi i \frac{m}{M}.}\chi_{E_l}\rangle \vert^2}.$$
				This for all $h\in \ell^2(E_l)$. Hence, $\{e^{2\pi i \frac{m}{M}.}\chi_{E_l}\}_{m\in \mathbb{N}_M, \, l\in \mathbb{N}_L}$ is a $M$-tight frame for $\ell^2(E)$.
			\end{enumerate}
		\end{proof}
		\begin{lemma}\label{lem4.4}
			Let $L,M\in \mathbb{N}$ and $E_0,\,E_1,\, \ldots,\,E_{L-1}\subset \mathbb{Z}$ be mutually disjoint. Denote $E=\displaystyle{\bigcup_{l\in \mathbb{N}_L}E_l}$. Then the following statements are equivalent:
			\begin{enumerate}
				\item $\{e^{2\pi i \frac{m}{M}.}\chi_{E_l}\}_{m\in \mathbb{N}_M, \, l\in \mathbb{N}_L}$ is a $M$-tight frame for $\ell^2(E)$. 
				\item $\{e^{2\pi i \frac{m}{M}.}\chi_{E_l}\}_{m\in \mathbb{N}_M, \, l\in \mathbb{N}_L}$ is complete in  $\ell^2(E)$.
				\item For all $l\in \mathbb{N}_L$, $E_l$ is $M\mathbb{Z}$-congruent to a subset of $\mathbb{N}_M$.
				\item For all $l\in \mathbb{N}_L$, $\displaystyle{\sum_{k\in \mathbb{Z}}\chi_{E_l}(.+kM)\leq 1}$ on $\mathbb{Z}$.
			\end{enumerate}
		\end{lemma}
		\begin{proof}
			It is a direct result of Lemma \ref{lemma6}, Lemma \ref{lem4.2} and Lemma \ref{lem4.3} together.
		\end{proof}

		We show in the following theorem that, for a $N\mathbb{Z}$-periodic set $\mathbb{S}$ in $\mathbb{Z}$, the admissibility of a  Gabor frame and the admissibility of a complete  Gabor system are equivalent and we provide a characterization based on the parameters $L,M,N$ and $R$.
		\begin{theorem}\label{prop6}
			The following statements are equivalent.
			\begin{enumerate}
				\item[(1)] There exists $g:=\{g_l\}_{l\in \mathbb{N}_L}\subset \ell^2(\mathbb{S},\mathbb{C}^R)$ such that $\mathcal{G}(g,L,M,N,R)$ is a Parseval frame for $\ell^2(\mathbb{S},\mathbb{C}^R)$.
				\item[(2)] There exists $g:=\{g_l\}_{l\in \mathbb{N}_L}\subset \ell^2(\mathbb{S},\mathbb{C}^R)$ such that $\mathcal{G}(g,L,M,N,R)$ is a frame for $\ell^2(\mathbb{S},\mathbb{C}^R)$.
				\item[(3)] There exists $g:=\{g_l\}_{l\in \mathbb{N}_L}\subset \ell^2(\mathbb{S},\mathbb{C}^R)$ such that $\mathcal{G}(g,L,M,N,R)$ is complete in $\ell^2(\mathbb{S},\mathbb{C}^R)$.
				\item[(4)] $R.card(\mathcal{K}_j)\leq q.L$ for all $j\in \mathbb{N}_{\frac{M}{q}}$.
			\end{enumerate}
		\end{theorem}
		\begin{proof}
			It is clear that $(1)\Longrightarrow (2)\Longrightarrow (3)$. Assume $(3)$, then by Proposition $\ref{prop4}$, we heve for all $j\in \mathbb{N}_{\frac{M}{q}}$ and a.e $\theta \in [0,1[$, $
			rank(Z_g(j,\theta))=R.card(\mathcal{K}_j)$. Since $Z_g(j,\theta)\in \mathcal{M}_{qL,pR}$, then $R.card(\mathcal{K}_j)\leq qL$. Hence $(3)\Longrightarrow (4)$. It remains to show that $(4)$ implies $(1)$.
			Fix $j\in \mathbb{N}_{\frac{M}{q}}$ and let $K$ be the maximal integer such that $K\displaystyle{\left\lfloor\frac{q}{R}\right\rfloor}\leq card(\mathcal{K}_j)$, where $\lfloor.\rfloor$ denotes the floor function. Define  for all $l\in \mathbb{N}_K$, $\mathcal{K}_j^l$; the set of $(l+1)$-th  $\displaystyle{\left\lfloor\frac{q}{R}\right\rfloor}$ elements of $\mathcal{K}_j$, $\mathcal{K}_j^K$; the set of the rest elements of $\mathcal{K}_j$ and for $l\in \mathbb{N}_L-\mathbb{N}_{K+1}$, we define $\mathcal{K}_j^l=\emptyset$. Then, for all $j\in \mathbb{N}_{\frac{M}{q}}$ and $l\in \mathbb{N}_L$ we have $R.card(\mathcal{K}_j^l)\leq q$. for $j\in \mathbb{N}_{\frac{M}{q}}$ and $l\in \mathbb{N}_L$, choose $\{s_{j,i}^r\}_{1\leq r\leq R,\, i\in \mathcal{K}_j^l}\subset \mathbb{N}_q$ such that $s_{j,i}^r\neq s_{j,i'}^{r'}$ if $(i,r)\neq (i',r')$ (this choice is guaranteed since $R.card(\mathcal{K}_j^l)\leq q$). Then for $j\in \mathbb{N}_{\frac{M}{q}}$, $l\in \mathbb{N}_L$ and $1\leq r\leq R$, we define:
			$$E_j^{l,r}=\left\lbrace
			\begin{array}{rcl}
				&\emptyset&\;\; \text{ if }\mathcal{K}_j^l=\emptyset,\\
				&\{j+k_{j,i}^l M+ s_{j,i}^rN:\; i\in \mathbb{N}_{card(\mathcal{K}_j^l)}\}& \;\; \text{ otherwise}.
			\end{array}
			\right.$$
			And thus $E^{l,r}:=\displaystyle{\bigcup_{j\in \mathbb{N}_{\frac{M}{q}}} E_j^{l,r}}$ and $g_{l,r}:=\chi_{E^{l,r}}$. Take for all $l\in \mathbb{N}_L$, $g_l:=(g_{l,1},g_{l,2},\ldots,g_{l,R}),$ and $g:=\{g_l\}_{l\in \mathbb{N}_L}$ and let's prove, using Corollary \ref{cor2}, that $\mathcal{G}(g,L,M,N,R)$ is a super tight frame with frame bound $M$. First of all, We prove that, for all $1\leq r\leq R$, $\mathcal{G}(\, \{g_{l,r}\}_{l\in \mathbb{N}_L},L,M,N)$ is a frame for $\ell^2(\mathbb{S})$. Fix $1\leq r\leq R$.  It suffices to prove that $E^{0,r},E^{1,r},\ldots, E^{L-1,r}$ are mutually disjoint, and for all $l\in \mathbb{N}_L$, $E^{l,r}$ is $M\mathbb{Z}$-congruent to a subset of $\mathbb{N}_M$, and $E^r:=\displaystyle{\bigcup_{l\in \mathbb{N}_L}E^{l,r}}$ is $N\mathbb{Z}$-congruent to $\mathbb{S}_N$.  Indeed, in this case, we have  $\ell^2(\mathbb{S})=\displaystyle{\bigoplus_{n\in \mathbb{Z}}\ell^2(E^r+nN)}$ and, by Lemma \ref{lem4.4},  $\{e^{2\pi i\frac{m}{M}.}\chi_{E^{l,r}}\}_{m\in \mathbb{N}_M,l\in \mathbb{N}_L}$ is a tight frame for $\ell^2(E^r)$ with frame bound $M$. Then,  for all $n\in \mathbb{Z}$, $\{e^{2\pi i\frac{m}{M}.}\chi_{E^{l,r}}(.-nN)\}_{m\in \mathbb{N}_M,l\in \mathbb{N}_L}$ is a tight frame for $\ell^2(E^r+nN)$ with frame bound $M$. Hence, by similar arguments used in the proof of Lemma \ref{lem4.3}, $\mathcal{G}(\,\{\chi_{E^{l,r}}\}_{l\in \mathbb{N}_L}\,,L,M,N):=\{e^{2\pi i\frac{m}{M}.}\chi_{E^{l,r}}(.-nN)\}_{n\in \mathbb{Z}, m\in \mathbb{N}_M,l\in \mathbb{N}_L}$ is a tight frame for $\displaystyle{\bigoplus_{n\in \mathbb{Z}}\ell^2(E^r+nN)}=\ell^2(\mathbb{S})$ with frame bound $M$.\\
			\begin{enumerate}
				\item[]$\rightarrow$ Let's show that for all $l\in \mathbb{N}_L$, $E^{l,r}$ is $M\mathbb{Z}$-congruent to a subset of $\mathbb{N}_M$. Let $l\in \mathbb{N}_L$. For this, it suffices to show that for all $j\in \mathbb{N}_{\frac{M}{q}},\, i\in \mathbb{N}_{card(\mathcal{K}_j^l)}$, we have: $$M|\, (j+k_{j,i}^l M-s_{j,i}^rN)-(j'+k_{j',i'}^lM-s_{j',i'}^rN)\;\Longrightarrow\; j=j' \text{ and }i=i'.$$
				Let $j,j'\in \mathbb{N}_{\frac{M}{q}}$ and $i,i'\in \mathbb{N}_{card(\mathcal{K}_j^l)}$ and suppose that $M|\, (j+k_{j,i}^l M-s_{j,i}^rN)-(j'+k_{j',i'}^lM-s_{j',i'}^rN)$. Then $M|\, j-j'+(k_{j,i}^l-k_{j',i'}^l)M-(s_{j,i}^r-s_{j',i'}^r)N$.\\
				Put $d=\displaystyle{\frac{M}{q}}$, then $M=dq$ and $N=dp$. Thus $dq|\, j-j'+(k_{j,i}^l-k_{j',i'}^l)dq-(s_{j,i}^r-s_{j',i'}^r)dp$, then $d|j-j'$. Hence, $j=j'$ since $j,j'\in \mathbb{N}_d$. On the other hand, we have $dq|(s_{j,i}^r-s_{j,i'}^r)dp$, then $q|(s_{j,i}^r-s_{j,i'}^r)p$, thus $q|s_{j,i}^r-s_{j,i'}^r$ since $p\wedge q=1$, hence $s_{j,i}^r=s_{j,i'}^r$ since $s_{j,i}^r,s_{j,i'}^r\in \mathbb{N}_q$. And then $i=i'$.\\
				Hence, for all $l\in \mathbb{N}_L$, $E^{l,r}$ is $M\mathbb{Z}$-congruent to a subset of $\mathbb{N}_M$.
				\item[]$\rightarrow$ Let's prove now that $E^r$ is $N\mathbb{Z}$-congruent to $\mathbb{S}_N$. We show first that $E^r$ is $N\mathbb{Z}$-congruent to a subset of $\mathbb{N}_N$. For this, let $(l,j,i),\, (l',j',i')\in \mathbb{N}_L\times \mathbb{N}_{\frac{M}{q}}\times \mathbb{N}_{card(\mathcal{K}_j)}$ and suppose that $N| (j-j')+(k_{j,i}^l-k_{j',i'}^{l'})M-(s_{j,i}^r-s_{j',i'}^r)N$. Put $d=\displaystyle{\frac{M}{q}}$, then $M=dq$ and $N=dp$. Thus $dp|\, j-j'+(k_{j,i}^l-k_{j',i'}^{l'})dq-(s_{j,i}^r-s_{j',i'}^r)dp$, then $d|j-j'$, hence $j=j'$ since  $j,j'\in \mathbb{N}_d$. On the other hand, we have $dp|(k_{j,i}^l-k_{j,i'}^{l'})dq$, then $p|k_{j,i}^l-k_{j,i'}^{l'}$, hence $k_{j,i}^l=k_{j,i'}^{l'}$
				since $k_{j,i}^l,k_{j,i'}^{l'}\in \mathbb{N}_p$. Then, $l=l'$ and $i=i'$ by definition of the elements $k_{j,i}^l$. Thus, $E^r$ is $N\mathbb{Z}$-congruent to a subset of $\mathbb{N}_N$. Observe that $E^r\subset \mathbb{S}$, then $E^r$ is $N\mathbb{Z}$-congruent to a subset of $\mathbb{S}_N$. By what above, we have, in particular, that for a fixed $r$, the $E_j^{l,r}$ are mutually disjoint (and also the $E^{l,r}$ are mutually disjoint). Then,  $$
				\begin{array}{rcl}
					card(E^r)&=&\displaystyle{\sum_{l\in \mathbb{N}_L}\sum_{j\in \mathbb{N}_{\frac{M}{q}}}card(\mathcal{K}_j^{l})}\\
					&=&\displaystyle{\sum_{j\in \mathbb{N}_{\frac{M}{q}}} card(\mathcal{K}_j)}\\
					&=&\displaystyle{\sum_{j\in \mathbb{N}_{\frac{M}{q}}}\sum_{n\in \mathbb{Z}}\chi_{\mathbb{S}_N}(j+\displaystyle{\frac{M}{q}}n)}\;\; \text{ Remark  } 2.5 \text{ in  \cite{4}}\\
					&=&\displaystyle{\sum_{j\in \mathbb{Z}}\chi_{\mathbb{S}_N}(j)}\\
					&=&card(\mathbb{S}_N).
				\end{array}$$
				Hence, $E^r$ is $N\mathbb{Z}$-congruent to $\mathbb{S}_N$.
			\end{enumerate}
			To complete the proof, by Corollary \ref{cor2}, it suffices to prove that $\mathcal{G}(\, \{g_{l,r}\}_{l\in \mathbb{N}_L},L,M,N)$, for $1\leq r\leq R$, are mutually strongly disjoint. Let $(s,k)\in \mathbb{N}_q	\times \mathbb{N}_p$, we have: $$\begin{array}{rcl}
				\left(Z_{g_{l,r}}(j,\theta)\right)_{s,k}&=&z_{pM}g_{l,r}(j+kM-sN,\theta)\\
				&=&\displaystyle{\sum_{n\in \mathbb{Z}}g_{l,r}(j+kM-sN+npM,\theta)e^{2\pi i n\theta}}\\
				&=&\displaystyle{\sum_{n\in \mathbb{Z}}\chi_{E^{l,r}}(j+kM-sN+npM,\theta)e^{2\pi i n\theta}}\\
				&=&\chi_{E^{l,r}}(j+kM-sN,\theta) \text{ by Theorem } 3.2 \text{ in \cite{4}}\\
			\end{array}$$
			Again by Theorem $3.2$ in \cite{4}, we have $\left(Z_{g_{l,r}}(j,\theta)\right)_{s,k}$ is nonzero if and only if $k\in \mathcal{K}_j^{l}$ and $s\in \{s_{j,i}^r\}_{i\in \mathbb{N}_{card(\mathcal{K}_j^l)}}$. Let $1\leq r\neq r'\leq R$. Since $\left(Z_{g_{l,r}}^*(j,\theta)Z_{g_{l,r'}}(j,\theta)\right)_{s,k}=\displaystyle{\sum_{n\in \mathbb{Z}}\left(Z_{g_{l,r}}^*(j,\theta)\right)_{n,s}\left(Z_{g_{l,r'}}(j,\theta)\right)_{n,k}}$ and that $s_{j,i}^r\neq s_{j,i'}^{r'}$ for all $i,i'\in \mathbb{N}_{card(\mathcal{K}_j^l)}$, then  $Z_{g_{l,r}}^*(j,\theta)Z_{g_{l,r'}}(j,\theta)=0$ for all $l\in \mathbb{N}_L$. Hence, by Proposition \ref{prop3.2}, $\mathcal{G}(\, \{g_{l,r}\}_{l\in \mathbb{N}_L},L,M,N)$, for $1\leq r\leq R$, are mutually strongly disjoint.
		\end{proof}

		\begin{remark}
			In the special case of $\mathbb{S}=\mathbb{Z}$, for all $j\in \mathbb{N}_{\frac{M}{q}}$, $\mathcal{K}(j)=\mathbb{N}_p$. Then the condition: $R.card(\mathcal{K}(j))= q.L$ for all $j\in \mathbb{N}_{\frac{M}{q}}$ is equivalent to $RN\leq LM$. Hence, in the case of $\mathbb{S}=\mathbb{Z}$, the following statements are equivalent:
			\begin{enumerate}
				\item[(1)] There exists $g:=\{g_l\}_{l\in \mathbb{N}_L}\subset \ell^2(\mathbb{Z},\mathbb{C}^R)$ such that $\mathcal{G}(g,L,M,N,R)$ is a super Gabor Parseval frame for $\ell^2(\mathbb{Z},\mathbb{C}^R)$.
				\item[(2)] There exists $g:=\{g_l\}_{l\in \mathbb{N}_L}\subset \ell^2(\mathbb{Z},\mathbb{C}^R)$ such that $\mathcal{G}(g,L,M,N,R)$ is a super Gabor frame for $\ell^2(\mathbb{Z},\mathbb{C}^R)$.
				\item[(3)] There exists $g:=\{g_l\}_{l\in \mathbb{N}_L}\subset \ell^2(\mathbb{Z},\mathbb{C}^R)$ such that $\mathcal{G}(g,L,M,N,R)$ is complete in $\ell^2(\mathbb{Z},\mathbb{C}^R)$.
				\item[(4)] $RN\leq LM$.
			\end{enumerate}
		\end{remark}

		The following result provides a necessary condition on the cardinality of \( \mathbb{S}_N \) for the existence of a Gabor frame, and characterizes when such a frame forms a Riesz basis in \( \ell^2(\mathbb{S}, \mathbb{C}^R) \).
		\begin{theorem}\label{prop7}
			Let $g:=\{g_l\}_{l\in \mathbb{N}_L}\subset \ell^2(\mathbb{S},\mathbb{C}^R)$,  then:
			\begin{enumerate}
				\item $\mathcal{G}(g,L,M,N,R)$ is a  Gabor frame for $\ell^2(\mathbb{S},\mathbb{C}^R)$ only when $R. card(\mathbb{S}_N)\leq LM.$
				\item Suppose that $\mathcal{G}(g,L,M,N,R)$ is a  Gabor frame for $\ell^2(\mathbb{S},\mathbb{C}^R)$. Then, the following statements are equivalent:
				\begin{enumerate}
					\item $\mathcal{G}(g,L,M,N,R)$ is a  Gabor Riesz Basis for $\ell^2(\mathbb{S},\mathbb{C}^R)$.
					\item $R.card(\mathbb{S}_N)=LM$. 
				\end{enumerate}
			\end{enumerate}
		\end{theorem}
		\begin{proof}\hspace{1cm}
			\begin{enumerate}
				\item Suppose that $\mathcal{G}(g,L,M,N,R)$ is a  Gabor frame for $\ell^2(\mathbb{S},\mathbb{C}^R)$. Then $\mathcal{G}(S^{\frac{-1}{2}} g,L,M,N,R)$ is a  Gabor Parseval frame for $\ell^2(\mathbb{S},\mathbb{C}^R)$. Then we can suppose that $\mathcal{G}(g,L,M,N)$ is a Parseval frame. Write for all $l\in \mathbb{N}_L$, $g_l=(g_{l,1},g_{l,2},\ldots, g_{l,R})$. Then for each $1\leq r\leq R$, $\mathcal{G}(\, \{g_{l,r}\}_{l\in \mathbb{N}_L},L,M,N)$ is a Parseval frame for $\ell^2(\mathbb{S})$. Then, by Lemma $\ref{lem2.18}$, for all $1\leq r\leq R$,  $\displaystyle{\sum_{l\in \mathbb{N}_L}\|g_{l,r}\|^2=\frac{card(\mathbb{S}_N)}{M}}$. Hence, $$\displaystyle{\sum_{r=1}^R\sum_{l\in \mathbb{N}_L}\|g_{l,r}\|^2=\frac{R.card(\mathbb{S}_N)}{M}}.$$ Since  $\displaystyle{\sum_{r=1}^R \|g_{l,r}\|^2=\|g_l\|^2}$, then:\begin{equation}
					\displaystyle{\sum_{l\in \mathbb{N}_L}\|g_l\|^2=\frac{R.card(\mathbb{S}_N)}{M}}.
				\end{equation} Since $E_{\frac{m}{M}}$ and $T_{nN}$ are both unitary operators of $\ell^2(\mathbb{S},\mathbb{C}^R)$ for all $m\in \mathbb{N}_M$ and $n\in \mathbb{Z}$, Then $1\geqslant \|E_{\frac{m}{M}}T_{nN}g_l\|=\|g_l\|$. Hence,  
				$$\displaystyle{\frac{R.card(\mathbb{S}_N)}{M}\leq L }.$$
				\item Assume that $\mathcal{G}(g,L,M,N,R)$ is a   Gabor frame for $\ell^2(\mathbb{S},\mathbb{C}^R)$. by the same argument in the $2$-th line of this proof, we can suppose that $\mathcal{G}(g,L,M,N,R)$ is a super Gabor  Parseval frame for $\ell^2(\mathbb{S},\mathbb{C}^R)$. Assume $(a)$, then by Lemma  \ref{lem2.17}, $\mathcal{G}(g,L,M,N,R)$ is a super Gabor orthonormal basis for $\ell^2(\mathbb{S},\mathbb{C}^R)$. By $(17)$, we have $L=\displaystyle{\frac{R.card(\mathbb{S}_N)}{M}}.$ Then $(a)\Longrightarrow (b)$. Conversely, assume $(b)$, then by $(17)$ again, we have $\displaystyle{\sum_{l\in \mathbb{N}_l}\| g_l\|^2=L}$. Since, by the same arguments above, $\|g_l\|\leq 1$ for all $l\in \mathbb{N}_L$, then for all $l\in \mathbb{N}_L$, $\|g_l\|=1$. Hence, Lemma \ref{lem2.17} completes the proof.
			\end{enumerate}
		\end{proof}

		The following lemma is very useful for establishing the theorem characterizing the admissibility of Gabor bases.
		\begin{lemma}\label{lem21}\hspace{1cm}
			\begin{enumerate}
				\item $R.card(\mathcal{K}(j))\leq q.L$ for all $j\in \mathbb{N}_{\frac{M}{q}}$ $\Longrightarrow$ $R. card(\mathbb{S}_N)\leq LM.$
				\item Assume that  $R.card(\mathcal{K}(j))\leq q.L$ for all $j\in \mathbb{N}_{\frac{M}{q}}$. Then, the following statements are equivalent.
				\begin{enumerate}
					\item $R. card(\mathbb{S}_N)= LM.$
					\item $R.card(\mathcal{K}(j))= q.L$ for all $j\in \mathbb{N}_{\frac{M}{q}}$
				\end{enumerate}
			\end{enumerate}
		\end{lemma}
		\begin{proof}\hspace{1cm}
			\begin{enumerate}
				\item Assume  that $R.card(\mathcal{K}_j)\leq \,qL$ for all $j\in \mathbb{N}_{\frac{M}{q}}$.  We have: 
				$$\begin{array}{rcl}
					R.card(\mathbb{S}_N)&=&\displaystyle{R.\sum_{j\in \mathbb{Z}}\chi_{\mathbb{S}_N}(j)}\\
					&=&\displaystyle{R.\sum_{j\in \mathbb{N}_{\frac{M}{q}}}\sum_{n\in \mathbb{Z}}\chi_{\mathbb{S}_N}(j+\displaystyle{\frac{M}{q}}n)}\\
					&=&\displaystyle{\sum_{j\in \mathbb{N}_{\frac{M}{q}}}R.card(\mathcal{K}_j)}\;\;\; \text{ by Remark } 2.5 \text{ in \cite{4}}\\
					&\leq& \displaystyle{\frac{M}{q}.qL=LM}.
				\end{array}$$
				\item Assume that $R.card(\mathcal{K}_j)\leq qL$ for all $j\in \mathbb{N}_{\frac{M}{q}}$.\\
				Assume that $R.card(\mathcal{K}_j)=qL$ for all $j\in \mathbb{N}_{\frac{M}{q}}$. Then by the proof of $1$, we have: $$\begin{array}{rcl}
					R.card(\mathbb{S}_N)&=&\displaystyle{\sum_{j\in \mathbb{N}_{\frac{M}{q}}}R.card(\mathcal{K}_j)}\\
					&=&\displaystyle{\frac{M}{q}.qL=LM}.
				\end{array}$$
				Conversely, assume that $R.card(\mathbb{S}_N)=\, LM$. Again by the proof of $1$, we have $\displaystyle{\sum_{j\in \mathbb{N}_{\frac{M}{q}}}R.card(\mathcal{K}_j)}=R.card(\mathbb{S}_N)=LM=\displaystyle{\frac{M}{q}.qL}$. Since 
				$R.card(\mathcal{K}_j)\leq \,qL$ for all $j\in \mathbb{N}_{\frac{M}{q}}$, then $R.card(\mathcal{K}_j)=\,qL$ for all $j\in \mathbb{N}_{\frac{M}{q}}$.
			\end{enumerate}
		\end{proof}
		
		Ultimately, we arrive at this theorem which characterizes the admissibility of  Gabor bases through the parameters $L,M,N$ and $R$ based on the previous results.
		\begin{theorem}\label{cor3}
			The following statements are equivalent.
			\begin{enumerate}
				\item[(1)] There exists $g:=\{g_l\}_{l\in \mathbb{N}_L}\subset \ell^2(\mathbb{S},\mathbb{C}^R)$ such that $\mathcal{G}(g,L,M,N,R)$ is a  Gabor  orthonormal basis for $\ell^2(\mathbb{S},\mathbb{C}^R)$.
				\item[(2)] There exists $g:=\{g_l\}_{l\in \mathbb{N}_L}\subset \ell^2(\mathbb{S},\mathbb{C}^R)$ such that $\mathcal{G}(g,L,M,N,R)$ is a  Gabor Riesz basis for $\ell^2(\mathbb{S},\mathbb{C}^R)$.
				\item[(3)] $R.card(\mathcal{K}_j)= q.L$ for all $j\in \mathbb{N}_{\frac{M}{q}}$.
			\end{enumerate}
		\end{theorem}
		\begin{proof}
			It is clear that $(1)\Longrightarrow (2)$. By Theorems \ref{prop6} and \ref{prop7} and Lemma \ref{lem21} together, we have clearly that $(2)\Longrightarrow (3)$. Assume $(3)$, then by proposition \ref{prop6}, there exists $g:=\{g_l\}_{l\in \mathbb{N}_L}\subset \ell^2(\mathbb{S},\mathbb{C}^R)$ such that $\mathcal{G}(g,L,M,N,R)$ is a Gabor Parseval frame for $\ell^2(\mathbb{S},\mathbb{C}^R)$. Lemma \ref{lem21} together with $(3)$ implies that $R.card(\mathbb{S}_N)=LM$. Hence, by Theorem \ref{prop7}, $\mathcal{G}(g,L,M,N,R)$ is a  Gabor Riesz basis for $\ell^2(\mathbb{S},\mathbb{C}^R)$. Thus, Lemma \ref{lem2.17} completes the proof.
		\end{proof}
		\begin{remark}
			In the special case of $\mathbb{S}=\mathbb{Z}$, for all $j\in \mathbb{N}_{\frac{M}{q}}$, $\mathcal{K}(j)=\mathbb{N}_p$. Then, the condition: $R.card(\mathcal{K}(j))= q.L$ for all $j\in \mathbb{N}_{\frac{M}{q}}$ is equivalent to $RN\leq LM$. Hence, in the case of $\mathbb{S}=\mathbb{Z}$, the following statements are equivalent:
			\begin{enumerate}
				\item[(1)] There exists $g:=\{g_l\}_{l\in \mathbb{N}_L}\subset \ell^2(\mathbb{Z},\mathbb{C}^R)$ such that $\mathcal{G}(g,L,M,N,R)$ is a  Gabor orthonormal basis for $\ell^2(\mathbb{Z},\mathbb{C}^R)$.
				\item[(2)] There exists $g:=\{g_l\}_{l\in \mathbb{N}_L}\subset \ell^2(\mathbb{Z},\mathbb{C}^R)$ such that $\mathcal{G}(g,L,M,N,R)$ is a  Gabor Riesz basis for $\ell^2(\mathbb{Z},\mathbb{C}^R)$.
				\item[(3)] $RN=LM$.
			\end{enumerate}
		\end{remark}
		
		We conclude with this example, which illustrates the importance of multi-window  Gabor frames as a good alternative to single-window  Gabor frames.
		\begin{example}
			In this example, we will use the notations already introduced in what above. 
			Let $R=2$, $M=4$, $N=6$ and $\mathbb{S}=\{3,\,5\}+6\mathbb{Z}=\{3,\,5,\,9,\,11\}+12\mathbb{Z}.$ We have then $p=3$ and  $q=2$, thus $\mathbb{N}_{\frac{M}{q}}=\mathbb{N}_2=\{0,1\}$. It is clear that $\mathcal{K}_0=\emptyset$ and $\mathcal{K}_1=\{1,\,2\}$. Then $\max\{card(\mathcal{K}_0),\;card(\mathcal{K}_1)\}=2$. Then, by Theorem \ref{prop6}, there is no  Gabor frame with a single window for $\ell^2(\mathbb{S},\mathbb{C}^2)$ but, by the same Theorem, the exisitence of a $L$-window  Gabor frame for $\ell^2(\mathbb{S},\mathbb{C}^2)$ is guaranteed for $L\geqslant 2$. Here we give an example of a two-window  Gabor frame. Define $g_0=(g_{0,1},g_{0,2}):=(\chi_{\{9\}},\chi_{\{3\}}), \, g_1=(g_{1,1},g_{1,2}):=(\chi_{\{5\}},\chi_{\{11\}})\in \ell^2(\mathbb{S},\mathbb{C}^2)$. By a simple computation, we get the following for a.e $\theta\in  [0,1[$:
			$Z_{g_{0,1}}(0,\theta)=\begin{pmatrix}
				0&0&0\\
				0&0&0
			\end{pmatrix},\;\; Z_{g_{0,2}}(0,\theta)=\begin{pmatrix}
				0&0&0\\
				0&0&0
			\end{pmatrix}, Z_{g_{1,1}}(0,\theta)=\begin{pmatrix}
				0&0&0\\
				0&0&0
			\end{pmatrix},\;\; Z_{g_{1,2}}(0,\theta)=\begin{pmatrix}
				0&0&0\\
				0&0&0
			\end{pmatrix}$, $Z_{g_{0,1}}(1,\theta)=\begin{pmatrix}
				0&0&1\\
				0&0&0
			\end{pmatrix},\;\; Z_{g_{0,2}}(1,\theta)=\begin{pmatrix}
				0&0&0\\
				0&0&1
			\end{pmatrix}, Z_{g_{1,1}}(1,\theta)=\begin{pmatrix}
				0&1&0\\
				0&0&0
			\end{pmatrix},$ and $Z_{g_{1,2}}(1,\theta)=\begin{pmatrix}
				0&0&0\\
				0&e^{2\pi i \theta}&0
			\end{pmatrix}$. Then, for a.e $\theta\in [0,1[$, we have: 
			$Z_{g_0}(0,\theta)=\begin{pmatrix}
				0&0&0&0&0&0\\
				0&0&0&0&0&0
			\end{pmatrix},\;\; Z_{g_1}(0,\theta)=\begin{pmatrix}
				0&0&0&0&0&0\\
				0&0&0&0&0&0
			\end{pmatrix},\;\;Z_{g_0}(1,\theta)=\begin{pmatrix}
				0&0&1&0&0&0\\
				0&0&0&0&0&1
			\end{pmatrix},$ and  $ Z_{g_1}(0,\theta)=\begin{pmatrix}
				0&1&0&0&0&0\\
				0&0&0&0&e^{2\pi i \theta}&0
			\end{pmatrix}.$ Let $x=(x_0,x_1,\ldots,x_5)^t\in \mathbb{C}^6$, we have for a.e $\theta\in [0,1[$, $\left\| Z_{g_0}(0,\theta)x\right\|^2=\left\| Z_{g_1}(0,\theta)x\right\|^2=0$, $\left\| Z_{g_0}(1,\theta)x\right\|^2=\vert x_2\vert^2+\vert x_5\vert^2$, $\left\| Z_{g_1}(1,\theta)x\right\|^2=\vert x_1\vert^2+\vert x_4\vert^2.$  Since $I_2\otimes \mathcal{K}(0)=0$, then $\left\| (I_2\otimes \mathcal{K}(0))x\right\|^2=0.$ We have also for a.e $\theta \in [0,1[$, $\left\| (I_2\otimes \mathcal{K}(1))x\right\|^2=\vert x_1\vert^2+\vert x_2\vert^2+\vert x_4\vert^2+\vert x_5\vert^2$. Then we have for all $j\in \{0,1\}$, a.e $\theta\in [0,1[$ and $x\in \mathbb{C}^6$,
			$$\left\langle Z_{g_0}^*(j,\theta)Z_{g_0}(j,\theta)x,x\right\rangle +\left\langle Z_{g_1}^*(j,\theta)Z_{g_1}(j,\theta)x,x\right\rangle = \left\langle (I_2\otimes \mathcal{K}(j))x,x\right\rangle.$$
			Hence, by Theorem \ref{prop5}, $\mathcal{G}(g,2,4,6,2)$ is a tight frame for $\ell^2(\mathbb{S},\mathbb{C}^2)$ with frame bound $4$, where $g:=\{g_0,g_1\}$. Observe that $R.card(\mathcal{K}_1)= q.L$ but $R.card(\mathcal{K}_0)\neq q.L$, then $\mathcal{G}(g,2,4,6,2)$  is not a super Riesz basis (Lemma \ref{lem21} together with Theorem \ref{prop7}). Moreover, Theorem \ref{cor3} shows that there is no  Gabor Riesz basis $\mathcal{G}(g,L,4,6,R)$ for $\ell^2(\mathbb{S},\mathbb{C}^R)$ for any $L,R\in \mathbb{N}$ since $card(\mathcal{K}_0)\neq card(\mathcal{K}_1)$.
		\end{example}
		\medskip
		\section*{Conclusion}
		
		To conclude, this paper provides a complete characterization of multi-window Gabor systems in the vector-valued discrete space \( \ell^2(\mathbb{S}, \mathbb{C}^R) \). In particular, we have established necessary and sufficient conditions under which a given family \( \mathcal{G}(g, L, M, N, R) \) forms a complete system or a frame in this space. Moreover, we have derived admissibility conditions on the periodic set \( \mathbb{S} \), expressed in terms of the structural parameters \( L, M, N, R \), ensuring the existence of complete Gabor systems, Parseval frames, or orthonormal bases.
		
		These results build upon and generalize the earlier work of Y.-Z. Li and Q.-F. Lian in \cite{3}, \cite{6}, and \cite{4}, whose contributions provided the foundation for our approach. Our framework extends their findings to a broader class of multi-window Gabor systems defined on periodically supported \(\ell^2\)-spaces with vector-valued sequences.
		
		\medskip
		
		More specifically, the following known results are recovered as particular cases of our general theorems:
		\begin{itemize}
			\item When \( R = 1 \) in Theorems~\ref{prop4}, \ref{prop5}, \ref{prop6}, and \ref{cor3}, we obtain Theorems 3.1, 4.1, 3.2, and 3.3, respectively, from \cite{3}.
			\item When \( L = 1 \), the same theorems reduce to Theorems 3.1, 3.2, 5.1, and Corollary 5.1, respectively, from \cite{6}.
			\item When \( R = 1 \) and \( L = 1 \), we recover Theorems 3.3, 3.4, 4.3, and 5.3, respectively, from \cite{4}.
		\end{itemize}
		\medskip
		This generalization not only unifies several previously known results but also opens new directions for research in time-frequency analysis on structured discrete spaces.
		
		\medskip
		\section*{Declarations}
		
		\medskip
		
		\noindent \textbf{Availablity of data and materials}\newline
		\noindent Not applicable.
		
		\medskip

		\noindent \textbf{Competing  interest}\newline
		\noindent The author declares that he has no competing interests.

		\medskip
		
		\noindent \textbf{Fundings}\newline
		\noindent  The author declares that there is no funding available for this article.

	\end{document}